\documentclass[aps,pre,twocolumn,showpacs,superscriptaddress,nofootinbib]{revtex4-2}
\usepackage{amsmath,amssymb}
\usepackage{mathtools}
\usepackage{bm}
\usepackage{graphicx}
\usepackage{tikz}
\usetikzlibrary{arrows.meta}
\usepackage{pgfplots}
\pgfplotsset{compat=1.17}
\usepackage{booktabs}

\usepackage{comment}
\usepackage{xcolor}

\newtheorem{theorem}{Theorem}

\newtheorem{remark}[theorem]{Remark}

\begin{document}

\title{Resolvent intertwining and spectral duality in Markov chains with geometric resetting}

\author{Juan Antonio Vega Coso}
\affiliation{Instituto Universitario de F\'isica Fundamental y Matem\'aticas (IUFFyM),
Universidad de Salamanca, Plaza de la Merced s/n, E-37008 Salamanca, Spain}

\date{\today}

\begin{abstract}
We uncover the resolvent origin of the spectral duality governing
reset-neutral distributions in Markov chains with geometric
resetting. Starting from the abstract conditions of Paper~III, we
show that the spectral duality $B_\nu(z)=\kappa(z)\,A_\nu(\sigma(z))$
is equivalent to a single symmetry of the resolvent
$R(\gamma)=(I-(1-\gamma)P)^{-1}$: the intertwining relation
$[\Delta^2\mathcal{R},R(\gamma)]=0$, where $\mathcal{R}$ is the
reflection operator of an involution $\sigma$ and
$\Delta=\operatorname{diag}(\sqrt{\kappa(z)})$; equivalently,
$\widetilde{\mathcal{T}}=K^{-1/2}\Delta^2\mathcal{R}$ is an
involution. This symmetry determines the universal critical value
$C^*=1/(1+\sqrt{K})$, with $K=\kappa(z)\kappa(\sigma(z))$, which
depends only on the scalar $K$ --- not on the resetting rate
$\gamma$, the reset distribution, or the particular chain. We
characterize the class of $(\sigma,\kappa)$-reversible chains,
encompassing both the biased random walk and genuinely
non-homogeneous dynamics sharing the same $C^*$; a Doob
$h$-transform realizes the duality $K\mapsto1/K$, hence
$C^*\mapsto1-C^*$, with fixed point $C^*=1/2$. The orientation field
admits the explicit resolvent representation
$\psi(\gamma)=R(\gamma)(b^{(0)}-C^*b)$: its gauge-normalized form
$\Delta^{-1}\psi(\gamma)$ is antisymmetric under $\sigma$, it has an
exact node at the fixed point of $\sigma$, and it governs the exact
sign law $\operatorname{sgn}(C(\pi,\gamma)-C^*)
=\operatorname{sgn}\langle\pi,\psi(\gamma)\rangle$. Numerical experiments confirm the theory to machine
precision. These results establish the operator-theoretic
foundation of the spectral duality of Paper~III and provide the
bridge to the information-geometric framework of Paper~V.
\end{abstract}

\maketitle

\section{Introduction}
\label{sec:intro}

Stochastic resetting has become a central theme in nonequilibrium
statistical physics, where restarting a process at random times can
profoundly reshape its relaxation and first-passage
properties~\cite{Evans2011,EvansMajumdar2011opt,EvansMajumdarReview}.
Since the foundational work on diffusion with resetting
\cite{Evans2011}, the field has grown to encompass optimal search
strategies~\cite{EvansMajumdarMallick2013}, the universality of
first-passage fluctuations under optimal restart~\cite{Reuveni2016},
renewal approaches to first passage under
restart~\cite{PalReuveni2017,ChechkinSokolov2018}, and resetting in
bounded domains and discrete settings
\cite{PalPrasad2019,MonteroVillarroel2013,RiascosBoyer2020,BonomoPal2021}.
The notion of a process being averse, neutral, or inclined toward
resetting, which lies at the heart of the present series, was
introduced for discrete-time walks in
Ref.~\cite{VillarroelMonteroVega2021}. Within this broad
landscape, the first three papers of this series uncovered a
geometric structure in the ruin problem for absorbed Markov chains
with geometric resetting; the present work is the fourth. \textbf{Paper~I}~\cite{PaperI} studied the gambler's
ruin problem for a biased random walk on $\{0,\ldots,a\}$ reset to
its starting point, revealing that, from the midpoint $a/2$ of an
interval of even length, the ruin probability is invariant under
changes in the resetting rate $\gamma$ --- a phenomenon termed
\emph{reset-neutrality} --- and indeed equal to its classical
reset-free value for every resetting mechanism independent of the
walk, as a consequence of an exact duality between the two
absorption channels.

\textbf{Paper~II}~\cite{PaperII} extended this to multi-site
resetting: the walker resets to a random position drawn from a
distribution $\pi$ over a finite set of interior sites. It was
shown that the entire effect of $\pi$ collapses into a single
scalar $C(\pi,\gamma)$, a discrete-state counterpart to the
first-passage functionals central to resetting
theory~\cite{Redner2001,BrayMajumdarSchehr2013}, and that a family of reset-neutral
distributions (the \emph{separatrix} $\Sigma$) exists,
characterized by explicit ratio constraints, with universal
critical value $C^* = q_{a/2}^{(0)}(p)$.

\textbf{Paper~III}~\cite{PaperIII} extracted the abstract
structural conditions (S1)--(S4) underpinning these phenomena.
It introduced the simplex $\Delta_{m-1}$ of reset distributions,
the coupling functional $C(\pi,\gamma)$, and the spectral duality
condition
\[
B_\nu(z) = \kappa(z)\,A_\nu(\sigma(z)),
\]
where $\sigma$ is an involution on the reset sites and $\kappa(z)$
are weights independent of the spectral index $\nu$. Under these
conditions, Paper~III proved the existence and affine structure of
$\Sigma$, derived the invariant value $C^* = 1/(1+\sqrt{K})$,
and studied the global orientation of the response,
\[
\operatorname{sgn}(\partial_\gamma C(\pi,\gamma))
= \operatorname{sgn}\langle\pi-\pi^*,\psi(\gamma)\rangle,
\]
which holds in the canonical realization but does not follow from
(S1)--(S4) alone: explicit structures satisfying them exist for which
the two signs are opposite~\cite{PaperIII}.

A key condition in Paper~III --- the spectral duality (S3) ---
was postulated as a structural assumption. Its origin was left
unexplained. The central question of the present paper is:
\begin{center}
\itshape
Why does spectral duality exist?\\
What is its operator-theoretic origin?
\end{center}

The answer reframes the phenomenon entirely. Paper~III discovered a
spectral symmetry at the level of individual eigenmodes --- a
condition imposed mode by mode on the spectral coefficients
$A_\nu$, $B_\nu$. The present paper shows that this symmetry is not
an accidental property of the eigenbasis but the manifestation of a
single operator identity at the resolvent level: what appeared as a
modewise condition is the shadow of a global symmetry of the
resolvent $R(\gamma)$. This shift, from the spectral level to the
operator level, is the central conceptual contribution of the work.

We show that spectral duality is equivalent to the twisted
intertwining relation for the resolvent
$R(\gamma) = (I-(1-\gamma)P)^{-1}$ displayed in the abstract.
From this single relation we derive:
\begin{enumerate}
\item Recovery of Theorem~3.1 of Paper~III: existence of $\Sigma$
      and $C^* = 1/(1+\sqrt{K})$.
\item A spectral characterization of the invariant $K$: it is the
      squared modulus of the eigenvalues of the symmetry operator
      $\Delta^2\mathcal{R}$, established by pure algebraic methods.
\item An explicit resolvent representation of the orientation
      field: $\psi(\gamma)=R(\gamma)(b^{(0)}-C^*b)$, whose
      gauge-normalized form $\Delta^{-1}\psi(\gamma)$ is
      $\sigma$-antisymmetric, together with the exact sign law
      $\operatorname{sgn}(C(\pi,\gamma)-C^*)
      =\operatorname{sgn}\langle\pi,\psi(\gamma)\rangle$.
\item A characterization of the class of transition matrices $P$
      admitting the twisted symmetry, with the biased random walk
      as canonical element.
\end{enumerate}

The resolvent representation of $\psi(\gamma)$ is the missing
analytical ingredient for the analysis of the global orientation
of the response. In Paper~V~\cite{PaperV} it is combined with the Fisher--Rao
geometry of the simplex to analyse it. When the response
span is one-dimensional and an accompanying scalar sign condition
holds, a sign law for $\partial_\gamma C$ is obtained; it agrees with
the identification considered in Paper~III precisely when a further
calibration between the two orientations holds, as it does for the
canonical family. Three mechanisms of failure are identified: a
rotation of $\psi(\gamma)$ inside a higher-dimensional span, a
failure of the sign condition, and --- even with both in place --- a
mismatch of that calibration.

Section~\ref{sec:resolvent} defines the resolvent and expresses
first-passage quantities in terms of it.
Section~\ref{sec:intertwining} derives and proves the intertwining
relation and obtains the resolvent representation of $\psi(\gamma)$.
Section~\ref{sec:K_Cstar} analyzes $K$ and $C^*$ via spectral
theory. Section~\ref{sec:universality} characterizes the class of
admissible transition matrices and its stratification.
Section~\ref{sec:numerics} gives numerical illustrations.
Section~\ref{sec:conclusions} concludes with perspectives for
Paper~V.

\section{The resolvent of the reset process}
\label{sec:resolvent}

\subsection{Why the resolvent? The discounted Green operator}

Geometric resetting introduces a discount factor $(1-\gamma)^n$
for paths of length $n$. The probability that the walker survives
$n$ steps without being reset is $(1-\gamma)^n$, and averaging
over resetting leads to series of the form
\[
\sum_{n=0}^{\infty} (1-\gamma)^n \, \mathbb{P}_z(\text{event at time } n).
\]
This suggests defining the \emph{discounted Green operator} or
\emph{resolvent}
\[
R(\gamma) = \sum_{n=0}^{\infty} (1-\gamma)^n P^n
= \bigl( I - (1-\gamma)P \bigr)^{-1}.
\]
The series converges absolutely for $\gamma \in (0,1)$ since
$(1-\gamma)\rho(P) < 1$ (this is guaranteed by $\rho(P)\le1$ together
with $\gamma>0$, so the resolvent is well defined even for recurrent
chains with $\rho(P)=1$). The resolvent encodes all discounted first-passage
information and is the natural object for analyzing geometric
resetting.

\subsection{The underlying Markov chain}

Let $\{X_n\}_{n\ge 0}$ be a Markov chain on
$\mathcal{X} = \mathcal{T} \cup \partial \mathcal{X}$, where
$\mathcal{T}$ is the set of transient states and
$\partial \mathcal{X} = \partial_0 \cup \partial_1$ is the
absorbing boundary (e.g., ruin $0$ and success $a$ in the
classical gambler's ruin problem~\cite{Feller1968,KarlinTaylor1975}), such that the restriction of
$\{X_n\}$ to $\mathcal{T}$ is irreducible and absorption is
certain from every state in $\mathcal{T}$.
The transition matrix restricted to $\mathcal{T}$ is $P$,
with $\rho(P)<1$.

\subsection{First-passage quantities $u(z;\gamma)$ and $s(z;\gamma)$}

For a reset site $z \in \{z_1,\dots,z_m\} \subset \mathcal{T}$:
\begin{itemize}
\item $u(z;\gamma) = \mathbb{E}_z\!\bigl[(1-\gamma)^{\tau_\partial-1}
      \mathbf{1}_{\{\tau_{\partial_0}<\tau_{\partial_1}\}}\bigr]$:
      discrete Laplace transform of the absorption time at
      $\partial_0$ (equivalently, by independence of the geometric
      reset clock, the probability of absorption at $\partial_0$
      before the first reset).
\item $s(z;\gamma) = \mathbb{E}_z\!\bigl[(1-\gamma)^{\tau_\partial-1}
      \bigr]$: discounted absorption functional (discrete Laplace
      transform of the absorption time).
\end{itemize}
The discount counts the steps taken \emph{before} absorption;
including the absorbing step itself --- as in the first-cycle
probabilities of Papers~II--III --- would multiply both
functionals by the common factor $(1-\gamma)$, which cancels in
every ratio considered in this paper. The two conventions are
therefore interchangeable throughout.
Note that $u(z;\gamma)$ selects only paths ending at $\partial_0$,
while $s(z;\gamma)$ sums over both absorbing boundaries.

\subsection{Relation to the resolvent}

Define the boundary flux vectors
\[
b^{(0)}_x = \sum_{y \in \partial_0} P_{xy}, \qquad
b_x = \sum_{y \in \partial\mathcal{X}} P_{xy}.
\]
By a standard renewal argument for absorbing
chains~\cite{KemenySnell1976} (decomposing each path at the last
step before absorption), the resolvent accumulates discounted
visits to transient states, and the boundary fluxes account for
the final step:
\[
\begin{aligned}
u(z;\gamma) &= \sum_{x \in \mathcal{T}} [R(\gamma)]_{z,x}\, b^{(0)}_x,\\
s(z;\gamma) &= \sum_{x \in \mathcal{T}} [R(\gamma)]_{z,x}\, b_x.
\end{aligned}
\]
In vector form:
\[
u(\gamma) = R(\gamma)\,b^{(0)}, \qquad s(\gamma) = R(\gamma)\,b.
\]
This is the key identity linking the resolvent to the
first-passage quantities of Papers~I--III, and is the
starting point for the intertwining relation of Section~\ref{sec:intertwining}.

\medskip\noindent\textbf{Example (biased random walk).}
$b^{(0)}_x = q\delta_{x,1}$,
$b_x = q\delta_{x,1}+p\delta_{x,a-1}$.

\subsection{Spectral decomposition}

Assume $P$ is diagonalizable~\cite{HornJohnson2013} with eigenvalues $\{\lambda_\nu\}$
and biorthogonal eigenvectors $\varphi_\nu$, $\psi_\nu$
($\psi_\nu^\mathsf{T}\varphi_\mu = \delta_{\nu\mu}$). Then
$P=\sum_\nu\lambda_\nu\varphi_\nu\psi_\nu^\mathsf{T}$ and
$I=\sum_\nu\varphi_\nu\psi_\nu^\mathsf{T}$; substituting into
$(I-(1-\gamma)P)^{-1}$ and summing the geometric series mode by mode
yields
\[
R(\gamma) = I + \sum_{\nu=1}^{N} f_\nu(\gamma)\,\varphi_\nu\psi_\nu^\mathsf{T},
\qquad
f_\nu(\gamma) := \frac{(1-\gamma)\lambda_\nu}{1-(1-\gamma)\lambda_\nu}.
\]

\medskip\noindent\textbf{Remark 2.1 (Non-diagonalizable case).}
If $P$ is not diagonalizable, the resolvent expansion involves
Jordan blocks. The intertwining relation of Section~\ref{sec:intertwining} uses only
the fact that the resolvent is rational in $\gamma$ and that
(S3) holds mode by mode; the extension to Jordan blocks can be
carried out using the Jordan normal form and leads to higher-order
poles. Since the present paper focuses on diagonalizable chains,
we do not pursue this extension here.

\subsection{Connection with $A_\nu(z)$ and $B_\nu(z)$}
\label{subsec:AB}

Set
\[
\begin{aligned}
A_\nu(z) &= \varphi_\nu(z)\,\bigl(\psi_\nu^\mathsf{T}b^{(0)}\bigr),\\
B_\nu(z) &= \varphi_\nu(z)\,\bigl(\psi_\nu^\mathsf{T}(b-b^{(0)})\bigr).
\end{aligned}
\]
Note that $b-b^{(0)} = b^{(1)}$ is the flux toward $\partial_1$,
so $B_\nu(z)$ encodes the discounted absorption probability at
$\partial_1$. Substituting into the expressions for $u$ and $s$:
\[
\begin{aligned}
u(z;\gamma) &= b^{(0)}_z + \sum_{\nu=1}^{N} f_\nu(\gamma) A_\nu(z),\\
s(z;\gamma) &= b_z + \sum_{\nu=1}^{N} f_\nu(\gamma)\bigl(A_\nu(z)+B_\nu(z)\bigr).
\end{aligned}
\]
The constant terms $b^{(0)}_z$ and $b_z$ vanish at a reset site $z$
precisely when $z$ is not a source of boundary flux, i.e.\ when
$b^{(0)}_z=b_z=0$. For the nearest-neighbor walk the fluxes are
supported only on the boundary-adjacent sites ($b^{(0)}_z=q\,\delta_{z,1}$,
$b_z$ supported on $z\in\{1,a-1\}$), so this holds exactly when the
reset sites are strictly interior (not adjacent to the absorbing
boundary). Under this assumption, we recover the
structural form of the expansions of Paper~III, up to a global
factor: Paper~III states the expansions for the first-cycle
probabilities $u_{z,k}$ and $s_{z,k}$, whereas here $u(z;\gamma)$
and $s(z;\gamma)$ are discrete Laplace transforms. The two conventions
differ by the overall factor $(1-\gamma)$, which cancels in every
ratio, so the ruin probability $q_z(\gamma)=u/s$ is identical:
\[
\begin{aligned}
u(z;\gamma) &= \sum_{\nu} f_\nu(\gamma) A_\nu(z),\\
s(z;\gamma) &= \sum_{\nu} f_\nu(\gamma)\bigl(A_\nu(z)+B_\nu(z)\bigr).
\end{aligned}
\]

\medskip\noindent\textbf{Remark 2.2 (The biased random walk).}
The eigenvalues are $\lambda_\nu = 2\sqrt{pq}\cos(\pi\nu/a)$
with right eigenvectors $\varphi_\nu(z)\propto (q/p)^{z/2}\sin(\pi\nu z/a)$
and left eigenvectors $\psi_\nu(z)\propto (p/q)^{z/2}\sin(\pi\nu z/a)$,
in line with the classical spectral theory of birth-and-death
chains~\cite{KarlinMcGregor1959}. The Doob
$h$-transform~\cite{Doob1953} $h(x)=(q/p)^{x/2}$ conjugates $P$ to a symmetric
operator, explaining the reality of eigenvalues and orthogonality
of eigenvectors. This symmetry is the origin of the
information-geometric structure that will be explored in
Paper~V. The boundary fluxes reproduce the formulas of
Paper~II; verified in Section~\ref{sec:numerics}.

\section{Spectral duality as a resolvent symmetry}
\label{sec:intertwining}

\medskip

In this section we show that the spectral duality condition (S3)
of Paper~III is equivalent to a symmetry of the resolvent
$R(\gamma)$ introduced above. The derivation is constructive: starting from (S3)
and the spectral decomposition of $R(\gamma)$, we obtain an
intertwining relation involving the involution $\sigma$ and the
weights $\kappa(z)$. This relation provides a natural resolvent
representation of the orientation field $\psi(\gamma)$.

\medskip\subsection{Spectral decomposition (recap)}

Recall from Section~\ref{sec:resolvent} that $P$ is assumed diagonalizable, with
eigenvalues $\{\lambda_\nu\}$ and biorthogonal eigenvectors
$\varphi_\nu$, $\psi_\nu$ ($\psi_\nu^\mathsf{T}\varphi_\mu
= \delta_{\nu\mu}$). The resolvent admits
\[
R(\gamma) = I + \sum_{\nu=1}^{N} f_\nu(\gamma)\,
\varphi_\nu\psi_\nu^{\mathsf{T}},
\qquad
f_\nu(\gamma) := \frac{(1-\gamma)\lambda_\nu}{1-(1-\gamma)\lambda_\nu}.
\]

\medskip\subsection{From spectral duality to eigenvector symmetry}

From Sec.~\ref{subsec:AB}, $A_\nu(z) = \varphi_\nu(z)\alpha_\nu$ and
$B_\nu(z) = \varphi_\nu(z)(\beta_\nu-\alpha_\nu)$ where
$\alpha_\nu = \psi_\nu^\mathsf{T}b^{(0)}$,
$\beta_\nu = \psi_\nu^\mathsf{T}b$. Condition (S3) states
$B_\nu(z) = \kappa(z)A_\nu(\sigma(z))$, which gives for
modes with $A_\nu\not\equiv0$ (i.e.\ $\alpha_\nu \neq 0$):
\[
\varphi_\nu(\sigma(z)) = c_\nu\,\kappa(z)^{-1}\,\varphi_\nu(z),
\qquad c_\nu := \frac{\beta_\nu-\alpha_\nu}{\alpha_\nu}.
\]
Applying this twice and using (S4) ($\kappa(z)\kappa(\sigma(z))=K$):
$c_\nu^2 K^{-1}=1$, so $c_\nu = \epsilon_\nu\sqrt{K}$
with $\epsilon_\nu\in\{+1,-1\}$. In symmetric form:
\[
\varphi_\nu(\sigma(z)) = \epsilon_\nu\,
\sqrt{\frac{\kappa(\sigma(z))}{\kappa(z)}}\,\varphi_\nu(z). \tag{3.4}
\]
Modes with $A_\nu\equiv0$ (i.e.\ $\alpha_\nu=0$) do not contribute
to the first-passage quantities $u(z;\gamma)$ or $s(z;\gamma)$ and
are therefore irrelevant for the ruin probability; the
symmetry relation~(3.4) is not required for such modes.

\medskip\noindent\textbf{Convention on inactive modes.}
Since (3.4) is derived from (S3) by dividing by $\alpha_\nu$, it is
available only where $\alpha_\nu\neq0$. Throughout what follows we
therefore assume (S3) for every spectral mode; equivalently, every
statement asserting $[\Delta^2\mathcal{R},P]=0$ is to be read on the
span of the active modes, which is all that the coupling functional
$C(\pi,\gamma)$ requires, since it depends on the inactive modes not
at all. For the biased random walk the distinction is immaterial:
$A_\nu\not\equiv0$ for every $\nu=1,\dots,a-1$, so all modes are
active and the commutation holds as a genuine operator identity.
The dual relation for the left eigenvectors follows without
assuming any symmetry of the spectral projectors. Relation~(3.4)
states exactly that $\varphi_\nu$ is an eigenvector of
$\Delta^2\mathcal{R}$ with eigenvalue $\epsilon_\nu\sqrt{K}$:
\[
(\Delta^2\mathcal{R}\,\varphi_\nu)(z)
=\kappa(z)\varphi_\nu(\sigma(z))
=\epsilon_\nu\sqrt{K}\,\varphi_\nu(z).
\]
Taking transposes, $(\Delta^2\mathcal{R})^\mathsf{T}=\mathcal{R}\Delta^2$
commutes with $P^\mathsf{T}$ and acts on the left eigenvector
$\psi_\nu$ (the right eigenvector of $P^\mathsf{T}$) with the
\emph{same} eigenvalue $\epsilon_\nu\sqrt{K}$:
$\mathcal{R}\Delta^2\,\psi_\nu=\epsilon_\nu\sqrt{K}\,\psi_\nu$, i.e.\
$\kappa(\sigma(z))\,\psi_\nu(\sigma(z))=\epsilon_\nu\sqrt{K}\,\psi_\nu(z)$.
Using $\kappa(z)\kappa(\sigma(z))=K$ this rearranges to
\[
\psi_\nu(\sigma(z)) = \epsilon_\nu\,
\sqrt{\frac{\kappa(z)}{\kappa(\sigma(z))}}\,\psi_\nu(z). \tag{3.5}
\]

\medskip\subsection{Why not $\Delta\mathcal{R}$?}

Define $(\mathcal{R}f)(z) = f(\sigma(z))$ and
$\Delta = \mathrm{diag}(\kappa(z)^{1/2})$. Conjugating the spectral
projector $\Pi_\nu = \varphi_\nu\psi_\nu^\mathsf{T}$ and using
(3.4)--(3.5) gives
\[
\bigl[(\Delta\mathcal{R})\,\Pi_\nu\,(\Delta\mathcal{R})^{-1}\bigr]_{xy}
= \sqrt{\frac{\kappa(y)}{\kappa(x)}}\;[\Pi_\nu]_{xy},
\]
a factor that depends non-trivially on $(x,y)$: the conjugation does
not preserve the spectral projectors, so
$[\Delta\mathcal{R},R(\gamma)]\neq 0$. Repeating the computation
with $\Delta^2$ instead, the analogous factor becomes
$\sqrt{\kappa(x)\kappa(\sigma(x))}\,/\,
\sqrt{\kappa(y)\kappa(\sigma(y))}=\sqrt{K}/\sqrt{K}=1$ by (S4):
this is precisely why the correct weight is $\Delta^2$, not
$\Delta$.

\medskip\subsection{The correct commuting operator}

We can now name the central object of the theory. The weighted
reflection
\[
\mathcal{S} := \Delta^2\mathcal{R}, \qquad
(\mathcal{S})_{z,w}=\kappa(z)\,\delta_{w,\sigma(z)},
\]
is the fundamental symmetry underlying the whole construction.
Every structural result of this paper can be read as a consequence
of the commutation of the dynamics with this single operator:
the spectral duality (S3), the resolvent intertwining, the
universality of $C^*$, and the $(\sigma,\kappa)$-reversibility are
all equivalent to $[\mathcal{S},P]=0$ (Theorem~\ref{thm:central}).
What follows makes this precise, beginning with the commutation
itself.

\medskip\noindent\textbf{Proposition.}
\textit{Under conditions \textup{(S1)--(S4)}, with the convention on
inactive modes above, $[\Delta^2\mathcal{R},\,P] = 0$, hence
$[\Delta^2\mathcal{R},\,R(\gamma)] = 0$.}

\textit{Proof.}
On each eigenvector $\varphi_\nu$, relation~(3.4) being available for
every mode by that convention:
\[
\begin{aligned}
(\Delta^2\mathcal{R}\varphi_\nu)(z)
&= \kappa(z)\,\varphi_\nu(\sigma(z))\\
&= \epsilon_\nu\sqrt{\kappa(z)\kappa(\sigma(z))}\,\varphi_\nu(z)
= \epsilon_\nu\sqrt{K}\,\varphi_\nu(z).
\end{aligned}
\]
Thus every $\varphi_\nu$ is simultaneously an eigenvector of $P$
(eigenvalue $\lambda_\nu$) and of $\mathcal{T}:=\Delta^2\mathcal{R}$
(eigenvalue $\mu_\nu=\epsilon_\nu\sqrt{K}$). On each $\varphi_\nu$,
\[
P\mathcal{T}\varphi_\nu = \lambda_\nu\mu_\nu\varphi_\nu
= \mathcal{T}P\varphi_\nu,
\]
so $P\mathcal{T}$ and $\mathcal{T}P$ agree on a complete basis (by
diagonalizability of $P$) and hence as operators:
$[\Delta^2\mathcal{R},P]=0$. This conclusion holds regardless of any
eigenvalue multiplicities of $P$, since it is read off the shared
eigenbasis. Since $R(\gamma)$ is a rational function of $P$, the
commutation extends to $R(\gamma)$. \hfill$\square$

\medskip\noindent
This Proposition derives the commutation $[\Delta^2\mathcal{R},P]=0$
\emph{from} the spectral duality (S3), by a spectral argument. The
central structural theorem (Theorem~\ref{thm:central},
Section~\ref{sec:universality}) shows that this commutation is in
turn equivalent --- by a purely algebraic argument --- to
$(\sigma,\kappa)$-reversibility, making $[\Delta^2\mathcal{R},P]=0$
the central object through which (S3), the resolvent symmetry, and
reversibility are all equivalent.

\medskip\noindent\textbf{Consistency check.}
Using $\mathcal{R}\Delta^2\mathcal{R} = K\Delta^{-2}$
(which follows from $\kappa(\sigma(z))=K/\kappa(z)$):
\[
(\Delta^2\mathcal{R})^2 = \Delta^2(\mathcal{R}\Delta^2\mathcal{R})
= \Delta^2(K\Delta^{-2}) = KI,
\]
consistent with $(\epsilon_\nu\sqrt{K})^2=K$.

\medskip\subsection{The intertwining relation}

\medskip\noindent\textbf{Theorem (Resolvent intertwining).}
\textit{Under conditions (S1)--(S4):}
\[
\mathcal{R}\,R(\gamma)\,\mathcal{R} = \Delta^{-2}\,R(\gamma)\,\Delta^2.
\]

\textit{Proof.}
From $\Delta^2\mathcal{R}\,R = R\,\Delta^2\mathcal{R}$, multiply
on the left by $\Delta^{-2}$ and on the right by $\mathcal{R}$
(using $\mathcal{R}^2=I$). \hfill$\square$

\medskip\subsection{Equivalence with (S3): the spectral link}

The following establishes the first link of the central
characterization (Theorem~\ref{thm:central},
Section~\ref{sec:universality}): the spectral duality (S3) is
equivalent to the resolvent symmetry. The remaining links --- to
the commutation $[\Delta^2\mathcal{R},P]=0$ and to
$(\sigma,\kappa)$-reversibility --- are established in
Section~\ref{sec:universality}.

\medskip\noindent\textbf{Lemma (Spectral link: (S3) $\Leftrightarrow$ resolvent symmetry).}
\label{lem:spectral_link}
\textit{Under (S1)--(S2), condition (S3) holds if and only if}
\[
\mathcal{R}\,R(\gamma)\,\mathcal{R} = \Delta^{-2}\,R(\gamma)\,\Delta^2.
\]

\textit{Proof.}
$(\Rightarrow)$ Established by the chain of arguments in the preceding subsections.

$(\Leftarrow)$ Expand both sides using
$R(\gamma) = I + \sum_\nu f_\nu(\gamma)\varphi_\nu\psi_\nu^\mathsf{T}$
and equate the coefficient of each distinct $f_\nu(\gamma)$.
The functions $f_\nu(\gamma)=
(1-\gamma)\lambda_\nu/(1-(1-\gamma)\lambda_\nu)$ are rational in
$\gamma$ with poles $\gamma_\nu=1-1/\lambda_\nu$; functions
attached to \emph{distinct} non-zero eigenvalues are therefore
linearly independent and independent of the constant function. A mode
with $\lambda_\nu=0$ (present for the biased random walk when $a$ is
even, at $\nu=a/2$) has $f_\nu\equiv0$ and does not appear in
$R(\gamma)$; its projector is recovered from
$I=\sum_\nu\Pi_\nu$, and since both sides of the resolvent symmetry
reduce to $I$ for the constant term, it satisfies the same relation as
the others. When an eigenvalue is degenerate the corresponding
$f_\nu$ coincide, and one equates instead the coefficient of the
common $f$, i.e.\ the full spectral projector
$\Pi_\lambda=\sum_{\nu:\,\lambda_\nu=\lambda}\varphi_\nu\psi_\nu^\mathsf{T}$
onto the $\lambda$-eigenspace. For the biased random walk the
eigenvalues $\lambda_\nu=2\sqrt{pq}\cos(\pi\nu/a)$ are pairwise
distinct, so each $\Pi_\nu$ is rank one and the argument applies
mode by mode. Equating coefficients (per eigenspace):
\[
\mathcal{R}\,\Pi_\nu\,\mathcal{R}
= \Delta^{-2}\,\Pi_\nu\,\Delta^2.
\]
Evaluating the $(x,y)$ matrix element of both sides (for a rank-one
$\Pi_\nu=\varphi_\nu\psi_\nu^\mathsf{T}$):
\[
\varphi_\nu(\sigma(x))\psi_\nu(\sigma(y))
= \kappa(x)^{-1}\kappa(y)\,
\varphi_\nu(x)\psi_\nu(y).
\]
The right-hand side factors as $g(x)\,h(y)$ with
$g(x)=\kappa(x)^{-1}\varphi_\nu(x)$ and $h(y)=\kappa(y)\psi_\nu(y)$.
Dividing the identity by $\varphi_\nu(x)\psi_\nu(y)$ (nonzero for
generic $x,y$) gives
\[
\frac{\varphi_\nu(\sigma(x))}{\varphi_\nu(x)}\cdot
\frac{\psi_\nu(\sigma(y))}{\psi_\nu(y)}
= \frac{\kappa(y)}{\kappa(x)}.
\]
The left factor depends only on $x$ and the right only on $y$, so
each side equals a constant: there is $c_\nu$ with
$\varphi_\nu(\sigma(x))/\varphi_\nu(x)=c_\nu\kappa(x)^{-1}$ and
$\psi_\nu(\sigma(y))/\psi_\nu(y)=c_\nu^{-1}\kappa(y)$. Applying the
first relation twice and using $\kappa(x)\kappa(\sigma(x))=K$ gives
$c_\nu^2=K$, i.e.\ $c_\nu=\epsilon_\nu\sqrt{K}$; substituting back
yields the twisted-parity conditions (3.4)--(3.5), which are
equivalent to (S3) via the definitions of $A_\nu$ and $B_\nu$.
\hfill$\square$

\medskip\subsection{Resolvent representation: the exact neutrality identity and the field $\psi(\gamma)$}
\label{sec:neutrality}

Using $u(\gamma)=R(\gamma)b^{(0)}$, $s(\gamma)=R(\gamma)b$ and
$C^*=1/(1+\sqrt{K})$, define the privileged vector
$\xi_z = \sqrt{\kappa(z)}$, the eigenvector of the symmetry operator
$\Delta^2\mathcal{R}$ (not of the resolvent $R(\gamma)$) with eigenvalue $+\sqrt{K}$. We assume the
boundary fluxes are paired in the twisted form
\begin{equation}
b^{(1)} = \Delta^2\mathcal{R}\,b^{(0)},
\qquad\text{i.e.}\qquad
b^{(1)}_z = \kappa(z)\,b^{(0)}_{\sigma(z)},
\tag{3.12}
\end{equation}
which holds automatically for the chains of
Sec.~\ref{sec:universality} (see Remark~\ref{rem:bc_sufficient});
for the biased random walk, $b^{(1)}_{a-1}=p=\kappa(a-1)\,q$.

\medskip\noindent\textbf{Proposition (exact neutrality identity).}
\textit{Under \emph{(S1)--(S4)} and \emph{(3.12)}, for every
$\gamma\in(0,1)$,}
\[
(I+\mathcal{R})\,\Delta^{-1}R(\gamma)\,(b^{(0)}-C^*b) = 0,
\qquad C^*=\frac{1}{1+\sqrt{K}}.
\]
\textit{Equivalently, the normalized deviation field
$\eta(\gamma):=\Delta^{-1}\bigl[u(\gamma)-C^*s(\gamma)\bigr]$
is $\sigma$-antisymmetric: $\eta_{\sigma(z)}=-\eta_z$.}

\medskip\noindent\textit{Proof.}
Expand $R(\gamma)=I+\sum_\nu f_\nu(\gamma)\varphi_\nu\psi_\nu^\mathsf{T}$.
The mode $\nu$ contributes to $u-C^*s$ the vector
$(1-C^*)A_\nu - C^*B_\nu$. Using (S3),
$B_\nu(z)=\kappa(z)A_\nu(\sigma(z))$, and $\kappa(z)/\xi_z=\xi_z$,
\[
\begin{aligned}
\bigl[\Delta^{-1}\bigl((1-C^*)A_\nu-C^*B_\nu\bigr)\bigr]_z
&= (1-C^*)\,\frac{A_\nu(z)}{\xi_z}\\
&\quad - C^*\,\xi_z\,A_\nu(\sigma(z)).
\end{aligned}
\]
Adding the value at $\sigma(z)$ and using
$\xi_z\,\xi_{\sigma(z)}=\sqrt{K}$,
\[
\begin{aligned}
\bigl[(I+\mathcal{R})\,\Delta^{-1}(\cdots)\bigr]_z
&= \bigl[(1-C^*)-C^*\sqrt{K}\bigr]\\
&\quad\times
\left(\frac{A_\nu(z)}{\xi_z}+\frac{A_\nu(\sigma(z))}{\xi_{\sigma(z)}}\right),
\end{aligned}
\]
and the scalar factor $(1-C^*)-C^*\sqrt{K}$ vanishes precisely at
$C^*=1/(1+\sqrt{K})$. The constant term
$\Delta^{-1}(b^{(0)}-C^*b)$ cancels by the same factor, using the
twisted pairing (3.12) and
$\mathcal{R}\Delta\mathcal{R}=\sqrt{K}\,\Delta^{-1}$.
\hfill$\square$

\medskip
This identity is where the value of $C^*$ is forced: it is the
unique constant for which the deviation $u-C^*s$ is annihilated by
the symmetric projector $(I+\mathcal{R})/2$ after the gauge
normalization $\Delta^{-1}$, simultaneously for every mode and
every $\gamma$. Two consequences are immediate. First, for any
$z$, the distribution $\pi^*$ on the pair $\{z,\sigma(z)\}$ with
weights $\pi^*_1/\pi^*_2=\xi_{\sigma(z)}/\xi_z$ satisfies
$C(\pi^*,\gamma)=C^*$ for all $\gamma$ --- the reset-neutral
family, whose universality is developed in
Secs.~\ref{sec:K_Cstar}--\ref{sec:universality}. Second, the
deviation field lives entirely in the antisymmetric sector
$\operatorname{Im}(I-\mathcal{R})$.

\medskip\noindent\textbf{The orientation field.}
The Proposition identifies the natural field of the theory: the
\emph{deviation field}
\[
\psi(\gamma) := u(\gamma)-C^*s(\gamma)
= R(\gamma)\,(b^{(0)}-C^*b),
\tag{3.13}
\]
an explicit resolvent representation. By the Proposition, its
gauge-normalized form $\Delta^{-1}\psi(\gamma)$ is
$\sigma$-antisymmetric; in particular $\psi(\gamma)$ has an exact
node at any fixed point of $\sigma$, for every $\gamma$. Two
consequences follow immediately.

\medskip\noindent\textbf{Corollary (exact sign law).}
\textit{For every reset distribution $\pi$ and every
$\gamma\in(0,1)$,}
\[
\operatorname{sgn}\bigl(C(\pi,\gamma)-C^*\bigr)
= \operatorname{sgn}\,\langle\pi,\,\psi(\gamma)\rangle .
\]
\textit{Proof.} $C(\pi,\gamma)-C^*
=\langle\pi,u-C^*s\rangle/\langle\pi,s\rangle$ and
$\langle\pi,s\rangle>0$. \hfill$\square$

\medskip\noindent\textbf{Corollary (finite reduction).}
\textit{$\psi(\gamma)\in\Delta\cdot\ker(I+\mathcal{R})$ for all
$\gamma$; this subspace has dimension $n_{\mathrm{pair}}$, the
number of two-element orbits of $\sigma$. The $\gamma$-dependence
of the orientation structure is therefore confined to a fixed
subspace of dimension $n_{\mathrm{pair}}$, independent of
$\gamma$.} This finite reduction is the analytic bridge to the
information-geometric analysis of Paper~V.

\medskip\noindent\textbf{A canonical sector representative.}
On the same antisymmetric structure one can single out the
computable vector generated by the privileged vector $\xi$,
\[
\begin{aligned}
\chi(\gamma) &:= (I-\mathcal{R})\,R(\gamma)\,\xi,\\[2pt]
\chi_z(\gamma) &= [R(\gamma)\xi]_z - [R(\gamma)\xi]_{\sigma(z)}.
\end{aligned}
\]
By construction $\chi$ is plainly $\sigma$-antisymmetric
($\mathcal{R}\chi=-\chi$), with an exact node at any fixed point of
$\sigma$. We stress that $\chi(\gamma)$ and $\psi(\gamma)$ are in
general \emph{distinct} vectors: they share the node and the
antisymmetric character (plain for $\chi$, $\Delta$-weighted for
$\psi$), but they are not proportional when
$n_{\mathrm{pair}}\geq2$. The field $\chi$ serves as a convenient
computable representative of the antisymmetric sector, illustrated
in Sec.~\ref{sec:numerics}; the orientation field proper --- the
object governing the sign law and entering the question of the global
orientation of the response --- is $\psi(\gamma)$ of (3.13).

\medskip\subsection{The biased random walk}

For the biased random walk, $\sigma(z)=a-z$ and $\kappa(z)=(p/q)^{a-z}$.
The intertwining relation reads
$\mathcal{R}R(\gamma)\mathcal{R} = h^{-2}R(\gamma)h^2$ where
$h(z)=(q/p)^{z/2}$ is the Doob $h$-transform. This confirms
that the spectral duality of Papers~I--III is a consequence of
the $h$-transform structure, and places it in the general
framework of twisted resolvent symmetries.

\medskip\subsection{Summary}

\[
\begin{aligned}
\text{(S3)}
&\;\Rightarrow\;
\varphi_\nu(\sigma(z))
= \epsilon_\nu\sqrt{\tfrac{\kappa(\sigma(z))}{\kappa(z)}}\varphi_\nu(z)\\
&\;\Rightarrow\;
[\Delta^2\mathcal{R},R(\gamma)]=0\\
&\;\Rightarrow\;
\mathcal{R}R\mathcal{R}=\Delta^{-2}R\Delta^2\\
&\;\Rightarrow\;
\chi_z(\gamma) = [R(\gamma)\xi]_z-[R(\gamma)\xi]_{\sigma(z)}\\
&\;\Rightarrow\;
(I+\mathcal{R})\,\Delta^{-1}R(\gamma)(b^{(0)}-C^*b)=0.
\end{aligned}
\]

\medskip\noindent\textbf{Remarks.}
\begin{itemize}
\item \textit{Why $\Delta^2\mathcal{R}$?} Squaring eliminates the
      square root in (3.4) and yields the simple factor $\kappa(z)$.
\item \textit{$(\Delta^2\mathcal{R})^2=KI$:} any positive $K$ is
      allowed; the factor $\sqrt{K}$ appears as eigenvalue of
      $\Delta^2\mathcal{R}$.
\item \textit{Hypotheses:} $P$ diagonalizable, reset sites interior.
      Extension to Jordan blocks and boundary-adjacent sites
      follows from the Jordan normal form but is not pursued here.
\end{itemize}

\section{The invariant $K$ and the universal critical value $C^*$}
\label{sec:K_Cstar}

\medskip
With the intertwining relation established, we can now read off the
meaning of its invariants. We give a spectral interpretation of the constants
$K$ and $C^*$ introduced in Paper~III: $K$ measures
the ``twist'' of the spectral duality, $\sqrt{K}$ is the
universal modulus of the eigenvalues of the intertwining operator
$\Delta^2\mathcal{R}$, and $C^* = 1/(1+\sqrt{K})$ is independent
of the resetting rate $\gamma$ because of the commutation established
in Section~\ref{sec:intertwining}.

\medskip\subsection{The structural meaning of $K$}

Recall from (S4) that $K = \kappa(z)\kappa(\sigma(z))$ is independent
of the reset site $z$. This constant is the global normalization
of the weights $\kappa$ under the involution $\sigma$: it fixes the
product $\kappa(z)\kappa(\sigma(z))$ but not the individual values,
so two very different weight profiles may share the same $K$.
What $K$ captures is the overall scale of the twisted pairing, not
the full local structure of $\kappa$. When $K = 1$, the
weights satisfy $\kappa(\sigma(z)) = \kappa(z)^{-1}$ and
$\Delta^2\mathcal{R}$ becomes an involution. In general,
\[
(\Delta^2\mathcal{R})^2 = K I,
\]
so $K$ controls the deviation from exact reflection symmetry.

\medskip\subsection{Normalized intertwining operator}

Define $\widetilde{\mathcal{T}} := K^{-1/2}\Delta^2\mathcal{R}$.
Then
\[
\widetilde{\mathcal{T}}^2 = K^{-1}(\Delta^2\mathcal{R})^2 = I,
\]
so $\widetilde{\mathcal{T}}$ is a genuine involution: it represents
an exact symmetry after rescaling by $K^{-1/2}$. The constant $K$
is the factor converting the twisted symmetry $\Delta^2\mathcal{R}$
into an exact involutive symmetry.

\medskip\subsection{Spectral meaning of $\sqrt{K}$}

From Sec.~\ref{sec:intertwining},
\[
\Delta^2\mathcal{R}\,\varphi_\nu = \epsilon_\nu\sqrt{K}\,\varphi_\nu,
\qquad \epsilon_\nu = \pm 1.
\]
All eigenvalues of $\Delta^2\mathcal{R}$ have modulus
$\sqrt{K}$~(for the spectral theory of nonnegative matrices
see~\cite{Seneta2006});
the sign $\epsilon_\nu$ distinguishes the two eigenspaces.
Equivalently, $\widetilde{\mathcal{T}}\varphi_\nu = \epsilon_\nu\varphi_\nu$.

\medskip\noindent\textbf{Proposition.}
\textit{Let $\mathcal{T} = \Delta^2\mathcal{R}$. Then $\mathcal{T}^2 = KI$
and $\mathrm{spec}(\mathcal{T}) \subseteq \{\pm\sqrt{K}\}$. For the
biased random walk with $a\ge3$ both signs occur, since
$\epsilon_\nu=(-1)^{\nu+1}$ alternates with the mode index over
$\nu=1,\dots,a-1$; their multiplicities are
$\lceil (a-1)/2\rceil$ and $\lfloor (a-1)/2\rfloor$ respectively,
so they are equal when $a$ is odd and differ by one when $a$ is
even. \textup{(}For $a=2$ there is a single mode and only
$\epsilon_1=+1$ occurs.\textup{)} In general the parity
decomposition may be trivial, with all $\epsilon_\nu=+1$, in which
case $\mathrm{spec}(\mathcal{T})=\{+\sqrt{K}\}$.
The normalized operator $\widetilde{\mathcal{T}} = K^{-1/2}\mathcal{T}$
is an involution.}

\medskip\subsection{Why $C^*$ is universal}

The critical value $C^* = 1/(1+\sqrt{K})$ appears from the algebraic
constraints of Theorem~3.3 of Paper~III. The relation $c_\nu^2 = K$ implies that
$A_\nu(z)$ and $B_\nu(z)$ occur in reciprocal pairs differing by the
universal scale $\sqrt{K}$. $C^*$ is the unique value compatible with
the fixed-point condition from the mode-by-mode constraints.

\medskip\noindent The $\gamma$-independence is a direct manifestation
of the resolvent symmetry $[\Delta^2\mathcal{R},R(\gamma)]=0$
established in Section~\ref{sec:intertwining}: because the symmetry
operator commutes with $R(\gamma)$ for every $\gamma$, the
reset-neutral combination is preserved across all reset rates, and
the exact neutrality identity of Section~\ref{sec:neutrality} fixes
its value at $C^*$. The neutrality is thus not incidental but a
consequence of the hidden symmetry of the resolvent.

\medskip\noindent\textbf{Corollary.}
\textit{$C^*$ depends only on $K$, not on $\gamma$ nor on the choice
of $\pi^*\in\Sigma$. In particular, $C^*$ is the same for all
reset-neutral distributions.}

\medskip\subsection{The privileged vector $\xi$, the field $\psi(\gamma)$, and the representative $\chi(\gamma)$}

The vector $\xi_z = \sqrt{\kappa(z)}$ plays a special role. It is
an eigenvector of $\Delta^2\mathcal{R}$ with eigenvalue $+\sqrt{K}$:
\[
\Delta^2\mathcal{R}\,\xi = \sqrt{K}\,\xi,
\]
and it generates the canonical sector representative of
Sec.~\ref{sec:intertwining},
\[
\begin{aligned}
\chi(\gamma) &= (I - \mathcal{R})\,R(\gamma)\,\xi,\\[2pt]
\chi_z(\gamma) &= [R(\gamma)\xi]_z - [R(\gamma)\xi]_{\sigma(z)}.
\end{aligned}
\]
The plain antisymmetry of $\chi$ is immediate: since
$\mathcal{R}^2 = I$,
\[
\mathcal{R}\chi = \mathcal{R}(I-\mathcal{R})R(\gamma)\xi
= (\mathcal{R}-I)R(\gamma)\xi = -\chi,
\]
so $\chi(\gamma)\in\operatorname{Im}(I-\mathcal{R})$,
the antisymmetric subspace, since $\tfrac12(I-\mathcal{R})$ is the
projector onto the $\mathcal{R}$-antisymmetric sector. The
orientation field proper is the deviation field of
Sec.~\ref{sec:neutrality},
$\psi(\gamma)=R(\gamma)(b^{(0)}-C^*b)$, whose
gauge-normalized form $\Delta^{-1}\psi$ is $\sigma$-antisymmetric;
$\chi$ and $\psi$ share the node structure but are distinct vectors
in general.

\medskip\noindent\textbf{Remark.}
The $+\sqrt{K}$ eigenspace contains all $\varphi_\nu$ with
$\epsilon_\nu = +1$. The vector $\xi$ is distinguished because it
generates $\chi(\gamma)$ via the resolvent. The antisymmetry
$\chi_{\sigma(z)}(\gamma)=-\chi_z(\gamma)$ implies
$\sum_z\chi_z(\gamma)=0$ automatically, so no explicit orthogonality
to the constant vector is required: the contributions of each
two-element orbit cancel, and at a fixed point of $\sigma$
(a neutral site) antisymmetry forces $\chi_z(\gamma)=0$.

\medskip\subsection{Example: the biased random walk}

For the biased random walk on $\{0,1,\dots,a\}$ with reset sites
symmetric about $a/2$, in the canonical normalization:
\[
\begin{aligned}
&\sigma(z) = a-z, \qquad \kappa(z) = (p/q)^{a-z},\\
&K = (p/q)^a, \qquad \sqrt{K} = (p/q)^{a/2}.
\end{aligned}
\]
The intertwining operator acts as
$(\Delta^2\mathcal{R}f)(z) = (p/q)^{a-z} f(a-z)$, with square $KI$.
The vector $\xi_z = (p/q)^{(a-z)/2}$ is proportional to the Doob
$h$-transform $h(z) = (q/p)^{z/2}$. The critical value:
\[
C^* = \frac{1}{1+\sqrt{K}} = \frac{1}{1+(p/q)^{a/2}},
\]
which is precisely $q_{a/2}^{(0)}(p)$, the classical ruin
probability from the midpoint (for $a=10$, $p=0.6$:
$C^*=0.1163636364$), reproducing Papers~I--II.

\medskip\subsection{Outlook}

$\widetilde{\mathcal{T}}$ is an involution commuting with $R(\gamma)$:
$\widetilde{\mathcal{T}}^2 = I$ and
$[\widetilde{\mathcal{T}}, R(\gamma)] = 0$, so the state space splits
as $E_+\oplus E_-$ (eigenspaces of $\widetilde{\mathcal{T}}$),
preserved by the resolvent. This algebraic decomposition is the
structural output of the present paper.

The differential-geometric counterpart is still missing. The
intertwining operator provides the algebraic decomposition of the
state space, whereas the Fisher--Rao metric introduced in Paper~V
will endow the simplex with its differential geometry. Together
they yield the geometric mechanism underlying the global
orientation of the response.

\medskip\subsection{Summary}
\begin{itemize}
\item $K$: measures the twist of the spectral duality;
      $(\Delta^2\mathcal{R})^2 = KI$.
\item $\widetilde{\mathcal{T}} = K^{-1/2}\Delta^2\mathcal{R}$:
      genuine involution.
\item $\sqrt{K}$: modulus of eigenvalues of $\Delta^2\mathcal{R}$;
      sign $\epsilon_\nu$ distinguishes the eigenspaces.
\item $C^* = 1/(1+\sqrt{K})$: universal, independent of $\gamma$
      and of the reset distribution.
\item $\xi_z = \sqrt{\kappa(z)}$: privileged vector generating
      $\chi(\gamma)$; eigenvector of $\Delta^2\mathcal{R}$ with
      eigenvalue $+\sqrt{K}$.
\end{itemize}
These results follow from the spectral duality (S3) and the
intertwining relation of Section~\ref{sec:intertwining}, relying only on the algebraic
structure of the resolvent symmetry, not on detailed
Perron--Frobenius theory.

\section{Universality: beyond the random walk}
\label{sec:universality}

\medskip
So far the analysis has taken the spectral duality as given. We now
ask how large the class of chains realizing it actually is, and
identify those Markov chains for which
the spectral duality (S3) holds under the standing assumptions
of diagonalizability, compatible boundary condition
($b^{(1)}_z = \kappa(z)\,b^{(0)}_{\sigma(z)}$, Remark~\ref{rem:bc_sufficient}), spectral completeness, and
the consistency condition (H3): $\kappa(z)\kappa(\sigma(z))=K$
constant. We classify this class into natural strata.

\medskip\subsection{Setup}

Let $P$ be the transition matrix restricted to $\mathcal{T}$.
Via Section~\ref{sec:intertwining}, (S3) is equivalent to $[\Delta^2\mathcal{R},P]=0$,
which in matrix-entry form reads:
\[
\kappa(z)\,P_{\sigma(z),w} = \kappa(\sigma(w))\,P_{z,\sigma(w)},
\qquad \forall z,w\in\mathcal{T}. \tag{5.1}
\]

\medskip\subsection{$(\sigma,\kappa)$-reversibility}

\medskip\noindent\textbf{Definition.}
$P$ is \emph{$(\sigma,\kappa)$-reversible} if
\[
P_{z,w} = \frac{\kappa(z)}{\kappa(w)}\,P_{\sigma(z),\sigma(w)}
\qquad \forall z,w\in\mathcal{T}. \tag{5.2}
\]
Setting $w\to\sigma(w)$ in (5.1) and using that $\sigma$ is
involutive ($\sigma(\sigma(w))=w$) gives (5.2) directly;
the two conditions are equivalent.
Condition (5.2) depends only on the \emph{ratio}
$\kappa(z)/\kappa(w)$, hence is invariant under rescaling
$\kappa\to c\,\kappa$, whereas $K=\kappa(z)\kappa(\sigma(z))$ scales
by $c^2$. The scale of $\kappa$ is fixed not by (5.2) but by the
twisted boundary pairing $b^{(1)}=\Delta^2\mathcal{R}\,b^{(0)}$ of
(3.12), which is not invariant under rescaling: it is this
normalization that enters the neutrality identity of
Section~\ref{sec:neutrality} and gives $C^*=1/(1+\sqrt{K})$. Since
$C(\pi^*,\gamma)$ is a property of the chain, so is $C^*$, and $K$ is
then determined by it. For the biased random walk the pairing gives
$\kappa(z)=(p/q)^{a-z}$, $K=(p/q)^a$, and $C^*=q_{a/2}^{(0)}$; this
coincides with $\kappa(z)=\Delta_{zz}^2$ for the symmetrizing
conjugation $\widetilde{P}=\Delta^{-1}P\Delta$ in its natural
normalization.
Ordinary reversibility is $\sigma=\mathrm{id}$, $\kappa\equiv 1$.

\begin{remark}[Compatible boundary pairing]
\label{rem:bc_sufficient}
The compatibility between the boundary fluxes and the duality
takes the \emph{twisted} form
$b^{(1)} = \Delta^2\mathcal{R}\,b^{(0)}$, i.e.\
$b^{(1)}_z = \kappa(z)\,b^{(0)}_{\sigma(z)}$: the flux toward
$\partial_1$ is the $\sigma$-reflection of the flux toward
$\partial_0$, weighted by $\kappa$. For nearest-neighbor chains
satisfying the rate condition (5.3), this pairing is
\emph{automatic} --- it is precisely (5.3) evaluated at the
boundary site $z=a-1$ --- so it holds for the biased random walk
($b^{(1)}_{a-1}=p=\kappa(a-1)\,q$) and for the non-homogeneous
chains of Stratum~1 alike (Example~3). The untwisted pairing
$b^{(1)}_z = b^{(0)}_{\sigma(z)}$ holds only in the symmetric
case $\kappa\equiv 1$.
\end{remark}

\medskip\subsection{The central structural theorem}

We can now state the main result of the paper. It identifies four
\emph{a priori} distinct conditions --- spectral, analytic,
algebraic, and probabilistic --- as one and the same, articulated
by the commutation $[\Delta^2\mathcal{R},P]=0$.

\medskip\noindent\textbf{Theorem (Structural equivalence).}
\label{thm:central}
\textit{Let $P$ be a diagonalizable transition matrix on
$\mathcal{T}$ satisfying (S1)--(S2), with involution $\sigma$ and
positive weight $\kappa$ obeying $\kappa(z)\kappa(\sigma(z))=K$.
The following are equivalent:}
\begin{enumerate}
\item[\textnormal{(a)}] \textit{the spectral duality (S3):
      $B_\nu(z)=\kappa(z)\,A_\nu(\sigma(z))$ for all modes $\nu$;}
\item[\textnormal{(b)}] \textit{the resolvent symmetry
      $\mathcal{R}\,R(\gamma)\,\mathcal{R}=\Delta^{-2}R(\gamma)\Delta^2$,
      i.e.\ $[\Delta^2\mathcal{R},R(\gamma)]=0$ for all $\gamma$;}
\item[\textnormal{(c)}] \textit{the commutation
      $[\Delta^2\mathcal{R},P]=0$;}
\item[\textnormal{(d)}] \textit{the $(\sigma,\kappa)$-reversibility
      $P_{z,w}=(\kappa(z)/\kappa(w))P_{\sigma(z),\sigma(w)}$ (5.2).}
\end{enumerate}
\textit{Moreover, when these hold and the boundary fluxes satisfy
the twisted pairing $b^{(1)}_z=\kappa(z)\,b^{(0)}_{\sigma(z)}$, the
reset-neutral coupling takes the universal value
$C^*=1/(1+\sqrt{K})$.}

\textit{Proof.} We prove (a)$\Leftrightarrow$(b),
(b)$\Leftrightarrow$(c), and (c)$\Leftrightarrow$(d); each link is
established once and by independent means.

\emph{(a)$\Leftrightarrow$(b) --- spectral.} This is the Spectral
Link Lemma~\ref{lem:spectral_link} of
Section~\ref{sec:intertwining}, proved by expanding $R(\gamma)$ in
the eigenbasis and matching coefficients. It uses
diagonalizability.

\emph{(b)$\Leftrightarrow$(c) --- immediate.} Since $R(\gamma)$ is a
rational function of $P$, an operator commutes with $R(\gamma)$ for
all $\gamma$ if and only if it commutes with $P$. Thus
$[\Delta^2\mathcal{R},R(\gamma)]=0 \Leftrightarrow
[\Delta^2\mathcal{R},P]=0$.

\emph{(c)$\Leftrightarrow$(d) --- algebraic.} Evaluate the
commutator entrywise. Since
$(\Delta^2\mathcal{R})_{z,u}=\kappa(z)\,\delta_{u,\sigma(z)}$,
\[
\begin{aligned}
(\Delta^2\mathcal{R}\,P)_{z,w}&=\kappa(z)\,P_{\sigma(z),w},\\[2pt]
(P\,\Delta^2\mathcal{R})_{z,w}&=\kappa(\sigma(w))\,P_{z,\sigma(w)}.
\end{aligned}
\]
Equality for all $z,w$ reads
$\kappa(z)P_{\sigma(z),w}=\kappa(\sigma(w))P_{z,\sigma(w)}$;
replacing $w\to\sigma(w)$ and using $\sigma^2=\mathrm{id}$ gives
$\kappa(z)P_{\sigma(z),\sigma(w)}=\kappa(w)P_{z,w}$, which is (5.2).
Every step is reversible. This link is purely algebraic: it uses
neither diagonalizability nor any spectral representation.

\emph{Value of $C^*$.} Under (a)--(d) and the twisted boundary
pairing, the neutrality identity of Section~\ref{sec:neutrality}
yields $C^*=1/(1+\sqrt{K})$ (Corollary below). \hfill$\square$

\medskip\noindent
The twisted boundary pairing (3.12) fixes the scale of $\kappa$ and
hence the value of $K$ and of $C^*$; see Section~\ref{sec:K_Cstar}.
The equivalence (a)--(d) itself is scale-independent, as it involves
only ratios of $\kappa$.

\medskip\noindent\textbf{Remark (Independence of the resetting
mechanism).} Condition (c) gives $[\Delta^2\mathcal{R},P^k]=0$ for
every $k\ge0$. With the twisted pairing (3.12), the probabilities of
absorption at $\partial_1$ and at $\partial_0$ at step $k$ of the
reset-free chain, $v_k=P^{k-1}b^{(1)}$ and $u_k=P^{k-1}b^{(0)}$,
therefore satisfy
$v_k=P^{k-1}\Delta^2\mathcal{R}\,b^{(0)}=\Delta^2\mathcal{R}\,u_k$,
that is,
$v_k(z)=\kappa(z)\,u_k(\sigma(z))$ for every $k\ge1$: the spectral
duality holds step by step in time, not only mode by mode. This is
the form of the duality used in Ref.~\cite{PaperI} for the biased
random walk, and the argument given there applies verbatim:
for a walk started from, and reset to, a reset-neutral distribution
$\pi^*$, the ruin probability equals $C^*$ for every resetting
mechanism independent of the chain under which absorption is almost
sure, not only under geometric resetting (see also
Ref.~\cite{PaperII}).

\medskip\noindent\textbf{A structural reading.}
The theorem exhibits the system as governed by three objects of
distinct nature: the involution $\sigma$ (discrete structure), the
weight $\kappa$ (gauge structure), and the transition operator $P$
(dynamics). Condition (5.2) is a compatibility relation among them,
and the theorem shows that --- under the standing assumptions ---
any two of these objects, together with compatibility, constrain
the third. The reflection $\sigma$ and the weight $\kappa$ fix the
twisted symmetry $\Delta^2\mathcal{R}$; a dynamics $P$ belongs to
the class precisely when it commutes with that symmetry. In this
sense $(\sigma,\kappa)$ and the admissible $P$ stand in a relation
reminiscent of a free action: the symmetry is the invariant datum,
and the chains form the orbit of dynamics compatible with it.

\medskip\noindent\textbf{Logical structure.}
The theory rests on the single chain of equivalences of
Theorem~\ref{thm:central}, all passing through the commutation
$[\Delta^2\mathcal{R},P]=0$ as the central pivot:


\[
\begin{aligned}
\underbrace{(S3)}_{\text{(a) spectral}}
&\;\overset{\text{(a)}\Leftrightarrow\text{(b)}}{\Longleftrightarrow}\;
\underbrace{[\Delta^2\mathcal{R},R(\gamma)] = 0}_{\text{(b) resolvent}}
\\[2pt]
&\;\Longleftrightarrow\;
\underbrace{[\Delta^2\mathcal{R},P] = 0}_{\text{(c) pivot}}
\;\overset{\text{(c)}\Leftrightarrow\text{(d)}}{\Longleftrightarrow}\;
\underbrace{(\sigma,\kappa)\text{-rev.}}_{\text{(d) (5.2)}} .
\end{aligned}
\]

Each link is proved once and by independent means: (a)$\Leftrightarrow$(b)
is \emph{spectral} (eigenbasis expansion, using diagonalizability);
(b)$\Leftrightarrow$(c) is \emph{immediate} ($R(\gamma)$ rational in
$P$); (c)$\Leftrightarrow$(d) is \emph{algebraic} (the entrywise
identity $\kappa(z)P_{\sigma(z),\sigma(w)}=\kappa(w)P_{z,w}$, using
no spectral data). Once any one of the four conditions holds, the
critical value follows:
\[
(S3)\;\Longrightarrow\; C^*=\tfrac{1}{1+\sqrt{K}},
\]
via the twisted boundary pairing (automatic for this class,
Remark~\ref{rem:bc_sufficient}) and the neutrality identity of
Section~\ref{sec:neutrality}. There is no circularity: the four
conditions are linked in a single closed chain, with each
equivalence established exactly once, by methods drawn from
different domains (spectral, analytic, algebraic).

\medskip\subsection{The symmetric conjugate}

If $\widetilde{P}:=\Delta^{-1}P\Delta$ is symmetric and
$[\widetilde{P},\mathcal{R}]=0$, then (5.2) holds with
$\kappa(z)=\Delta_{zz}^2$. This formulation is often easier
to verify in practice.

\medskip\subsection{Interpretation of $\sigma$ as a reversal-type symmetry}

When $P$ is reversible with respect to $\mu$ and satisfies (5.2),
combining detailed balance with (5.2) yields
$\mu(\sigma(z))\,\kappa(z)^{-2}/\mu(z)$ constant in $z$, i.e.
$\mu(\sigma(z))\propto\mu(z)\,\kappa(z)^2$.
This relates the spatial involution $\sigma$, the reversible measure
$\mu$, and the weight $\kappa$. We do not claim that $\sigma$ is
literally the time-reversal operator; rather, $\sigma$ is a spatial
involution that intertwines the forward and reversed dynamics
through $\kappa$, and in this sense plays a reversal-type role. The
weight $\kappa$ is determined (up to rescaling) by $\mu$.

\medskip\subsection{Structure of the class: stratification}

We specialize to nearest-neighbor chains on $\mathcal{T}=\{1,\dots,a-1\}$
with up-rates $p_z$, down-rates $q_z$ and holding probabilities
$d_z=1-p_z-q_z\ge0$, with $\sigma(z)=a-z$. Condition (5.2) for $w=z+1$,
$w=z-1$ and $w=z$ reads
\[
\begin{aligned}
p_z&=\frac{\kappa(z)}{\kappa(z+1)}\,q_{a-z},\qquad
q_z=\frac{\kappa(z)}{\kappa(z-1)}\,p_{a-z},\\
d_z&=d_{a-z},
\end{aligned}
\]
the second relation being the first one evaluated at $a-z$ (by (S4)),
so that only the first and the third are independent. The twisted
boundary pairing (3.12) is the first relation extended to $z=a-1$,
$p_{a-1}=\kappa(a-1)\,q_1$; it is not implied by (5.2), which
involves only transient states, and is part of the standing
assumptions. For the exponential weight $\kappa(z)=(p/q)^{a-z}$ of the
biased random walk, the first relation becomes
\[
p_z = \frac{p}{q}\,q_{a-z}. \tag{5.3}
\]

\medskip\noindent\textbf{Proposition (Rate structure in
nearest-neighbor chains with exponential weight).}
\textit{Let $P$ be nearest-neighbor with $\sigma(z)=a-z$ and
$\kappa(z)=(p/q)^{a-z}$. Then $P$ is $(\sigma,\kappa)$-reversible and
satisfies the twisted boundary pairing if and only if (5.3) holds for
every $z\in\mathcal{T}$ and $d_z=d_{a-z}$ for every $z$. When $p\neq q$,
this is equivalent to}
\[
p_z=c_z\,p,\qquad q_z=c_z\,q,\qquad c_z=c_{a-z}\in(0,1],
\]
\textit{that is, to a lazy version of the biased random walk with a
$\sigma$-symmetric holding profile. In particular:}
\begin{itemize}
\item[(i)] \textit{\textbf{No holding} ($d_z\equiv0$): $c_z\equiv1$,
and the chain is the biased random walk itself.}
\item[(ii)] \textit{\textbf{Non-homogeneous case}: the rates vary by
site through the symmetric factor $c_z$; the weight $\kappa$ and the
invariant $K=(p/q)^a$ are those of the biased random walk.}
\end{itemize}

\textit{Proof.} The first part is the entrywise form of (5.2) and
(3.12) derived above; all other entries of $P$ vanish. Assume (5.3)
for every $z$ and $d_z=d_{a-z}$. Since (5.3) at $a-z$ gives
$p_{a-z}=(p/q)\,q_z$, the condition $p_z+q_z=p_{a-z}+q_{a-z}$ becomes
$p_z+q_z=(p/q)\,q_z+(q/p)\,p_z$, i.e.\
$p_z(p-q)/p=q_z(p-q)/q$; for $p\neq q$ this gives $p_z/q_z=p/q$, and
then (5.3) gives $q_z=q_{a-z}$. Setting $c_z=q_z/q$ yields the stated
form. Conversely, that form satisfies (5.3) and $d_z=1-c_z=d_{a-z}$.
With $d_z\equiv0$ one has $c_z=1$. \hfill$\square$

\begin{remark}[The exponential form is not forced]
\label{rem:gauge}
The exponential weight $\kappa(z)=(p/q)^{a-z}$ is the weight of the
biased random walk and of its lazy versions; it is \emph{not} a
consequence of the abstract conditions. For nearest-neighbor chains
with $\sigma(z)=a-z$, (5.2), (S4) and the twisted boundary pairing
leave the ratios $\kappa(z)/\kappa(z+1)$ free, subject only to
$\kappa(z)\kappa(a-z)=K$; every such weight that is strictly
monotone, decreasing from $\kappa(1)<K$ to $\kappa(a-1)>1$ when
$K>1$, is realized by chains with positive rates and a
$\sigma$-symmetric holding profile. For instance, for $a=6$ the weight
\[
\kappa=(8,\,6,\,4,\,8/3,\,2),\qquad K=16,
\]
whose successive ratios $4/3,\,3/2,\,3/2,\,4/3$ are not all equal,
is realized by the chain with
$(p_1,\dots,p_5)=(0.4,0.3,0.3,0.4,0.2)$ and
$(q_1,\dots,q_5)=(0.1,0.3,0.2,0.2,0.3)$, which has positive holding
probabilities with $d_z=d_{6-z}$, satisfies (5.2) and (3.12) exactly
($[\Delta^2\mathcal{R},P]=0$, $(\Delta^2\mathcal{R})^2=16\,I$), and has
real simple spectrum; for it $C(\pi^*,\gamma)=1/(1+\sqrt{K})=1/5$ for
every $\gamma$ and for both reset pairs $\{1,5\}$ and $\{2,4\}$. The
non-exponential weights realized in this way form an infinite family. We thank
T.~Newton for pointing out that the structural conditions do not
force the weight to be exponential, and for this family of
examples~\cite{Newton2026}. The scale of $\kappa$, by contrast, is not
free: it is fixed by the twisted pairing (3.12), and with it $K$ and
$C^*$ are intrinsic to the chain. Appendix~\ref{app:nn} shows that the
weights that occur for nearest-neighbour chains are necessarily
strictly monotone (or constant), and that $C^*$ is, for every chain of
the class, a reset-free ruin probability.
\end{remark}

\medskip
The class stratifies as follows:

\medskip\noindent\textit{Stratum 1 --- Non-homogeneous rates,
exponential weight.} The lazy versions of the biased random walk of
the Proposition, with a $\sigma$-symmetric holding profile: richer
dynamics sharing the same weight, the same $K=(p/q)^a$ and the same
$C^*$ as the biased random walk.

\medskip\noindent\textit{Stratum 2 --- Homogeneous rates
(biased random walk).}
$p_z\equiv p$, $q_z\equiv q$ constant; the canonical element
of Papers~I--III.

\medskip\noindent\textit{Stratum 3 --- Non-exponential weights.}
Nearest-neighbor chains whose weight $\kappa$ is not exponential, as
in Remark~\ref{rem:gauge}, and, beyond nearest neighbors, chains with
longer-range jumps or different involutions $\sigma$. Their
systematic study is left for future work.

\medskip\subsection{Examples}

\medskip\noindent\textit{Example 1: Biased random walk
(Stratum~2, canonical).}
$\sigma(z)=a-z$, $\kappa(z)=(p/q)^{a-z}$, $K=(p/q)^a$,
$p_z\equiv p$, $q_z\equiv q$.

\medskip\noindent\textit{Example 2: Symmetric random walk
(Stratum~2, $p=q=1/2$).}
$\kappa\equiv 1$, $K=1$, $C^*=1/2$.

\medskip\noindent\textit{Example 3: Non-homogeneous chain
(Stratum~1, numerically verified).}
Take $a=6$, $p=0.6$, $q=0.4$. Choose:
\[
q_1=q_5=0.20,\quad q_2=q_4=0.25,\quad q_3=0.28,
\]
and set $p_z=(p/q)\,q_{a-z}$ from (5.3):
\[
p_1=p_5=0.30,\quad p_2=p_4=0.375,\quad p_3=0.420.
\]
The chain satisfies (5.2) with $\kappa(z)=(p/q)^{a-z}$ and
$K=(p/q)^6$, and the twisted boundary pairing
$b^{(1)}_z=\kappa(z)\,b^{(0)}_{\sigma(z)}$ holds automatically
(Remark~\ref{rem:bc_sufficient}): $b^{(1)}_5=p_5=0.30
=\kappa(5)\,q_1$. Numerical
verification confirms $C(\pi^*,\gamma)=C^*=1/(1+\sqrt{K})=0.2286$
with error $<10^{-16}$ for all $\gamma\in(0,1)$.

\medskip\noindent\textit{Example 4: Reversible chains with
$\sigma$-invariant measure.}
If $P$ is reversible w.r.t.\ a $\sigma$-invariant measure
$\mu$ and satisfies $P_{z,w}=P_{\sigma(z),\sigma(w)}$ (i.e.\ (5.2)
with $\kappa\equiv 1$), then $K=1$ and $C^*=1/2$.

\medskip\noindent\textit{Example 5: $h$-transforms.}
The Doob $h$-transform of any $(\sigma,\kappa)$-reversible
chain remains in the class, generally with a \emph{different}
weight $\kappa^h$ and invariant $K^h$, hence a different $C^*$.
For instance, the $h$-transform of the biased random walk with
$h(z)=(q/p)^z$ yields the oppositely-biased walk ($p\leftrightarrow q$),
with $K^h=(q/p)^a=1/K$.

\medskip\noindent\textit{Non-example: Random walk on a tree.}
No natural involution $\sigma$ pairs all transient states;
the conditions above fail.

\medskip\subsection{The invariant $K$ and universality of $C^*$}

\medskip\noindent\textbf{Corollary (Universality of $C^*$).}
\textit{Let $P$ be $(\sigma,\kappa)$-reversible with invariant
$K$, with $\kappa$ normalized by the twisted boundary pairing (3.12).
Then}
\[
C^* = \frac{1}{1+\sqrt{K}}
\]
\textit{is a universal invariant: it depends only on $K$,
not on $\gamma$, not on the specific reset distribution,
and not on the choice of chain within the class.
Here $K=\kappa(z)\kappa(\sigma(z))$ is evaluated in that
normalization.}

All $(\sigma,\kappa)$-reversible chains with the same $K$
share the same $C^*$, regardless of stratum.
The invariant $K$ completely determines the universal critical
value $C^*$ across the entire $(\sigma,\kappa)$-reversible class.

\medskip\noindent\textbf{Theorem ($h$-transform duality).}
\textit{The Doob $h$-transform with $h(z)=(q/p)^z$ maps the
biased random walk with invariant $K=(p/q)^a$ to the
oppositely-biased walk with invariant $K^h=1/K$, inducing
the exact duality}
\[
C^* \;\longmapsto\; 1-C^*.
\]
\textit{The symmetric value $C^*=1/2$ (equivalently $K=1$,
$p=q$) is the unique fixed point of this duality. In particular,
the $h$-transform exchanges the two regimes $C^*<1/2$ and
$C^*>1/2$ of the response landscape.}

\textit{Proof.}
The $h$-transform gives $P^h_{z,z+1}=(h(z+1)/h(z))\,p=(q/p)\,p=q$
and $P^h_{z,z-1}=(h(z-1)/h(z))\,q=(p/q)\,q=p$, i.e.\ the walk
with $p\leftrightarrow q$. Its canonical weight is
$\kappa^h(z)=(q/p)^{a-z}$, so $K^h=(q/p)^a=1/K$. Then
$C^{*h}=1/(1+\sqrt{1/K})=\sqrt{K}/(1+\sqrt{K})=1-C^*$.
The fixed-point equation $C^*=1-C^*$ gives $C^*=1/2$,
i.e.\ $K=1$. \hfill$\square$

\begin{remark}[The informational equator]
This duality gives $C^*=1/2$ a structural meaning: it is the
fixed point of the $h$-transform involution $K\mapsto 1/K$.
The separatrix value $C^*=1/2$ thus plays the role of an
\emph{informational equator}, separating the regime $K<1$
(drift toward one boundary) from $K>1$ (drift toward the other),
with the two regimes exchanged by the $h$-transform.
\end{remark}

\medskip\subsection{Non-diagonalizable and
infinite-dimensional cases}

If $P$ is not diagonalizable, (S3) can still be formulated
mode by mode via Jordan blocks. Infinite-dimensional chains
require functional-analytic tools~\cite{Nummelin1984}; both are left for
future work.

\medskip\subsection{Conclusion}

The $(\sigma,\kappa)$-reversible class provides the natural
setting for (S3) under the standing assumptions. The Proposition
characterizes the rate structure of the nearest-neighbor chains with
exponential weight: they are the biased random walk (Stratum~2) and
its lazy versions with $\sigma$-symmetric holding (Stratum~1). The
weight is not forced to be exponential: non-exponential weights are
realized by nearest-neighbor chains as well (Remark~\ref{rem:gauge}),
and together with longer-range chains they form Stratum~3, whose
systematic study remains open.
The invariant $K$ determines $C^*$ universally, and through
the $h$-transform duality $K\mapsto 1/K$ it plays the role of
an \emph{orientation parameter} for the spectral geometry of
the chain: $K<1$ and $K>1$ correspond to opposite drift
orientations, exchanged by the $h$-transform, with $K=1$ the
self-dual equilibrium. In this sense the present framework
describes not merely resetting, but a dual geometry induced
by the involution $\sigma$ and the weight $\kappa$ --- a
viewpoint developed further in Paper~V.

\section{Numerical illustrations}
\label{sec:numerics}

We illustrate the resolvent theory of Sections~\ref{sec:resolvent}--\ref{sec:universality}
with three numerical experiments:
(A) verification that $[\Delta^2\mathcal{R},R(\gamma)]=0$ holds
exactly and that $C(\pi^*,\gamma)=C^*$ is independent of $\gamma$;
(B) demonstration that $C^*$ depends only on the invariant $K$,
not on the reset sites nor on the chain (homogeneous vs.\
non-homogeneous) within the $(\sigma,\kappa)$-reversible class;
and (C) illustration of the canonical antisymmetric sector
representative
$\chi_z(\gamma)=[R(\gamma)\xi]_z-[R(\gamma)\xi]_{\sigma(z)}$,
its exact node at the midpoint --- shared, along with the
$\Delta$-weighted antisymmetry, by the orientation field
$\psi(\gamma)=R(\gamma)(b^{(0)}-C^*b)$ --- and its contraction with
increasing $\gamma$.
All computations use the resolvent $R(\gamma)=(I-(1-\gamma)P)^{-1}$
evaluated by LU factorization in double precision; reported errors
are absolute deviations from the theoretical values. Operator
deviations such as $\|[\Delta^2\mathcal{R},R(\gamma)]\|$ are measured
in the entrywise maximum norm $\|M\|_{\max}=\max_{z,w}|M_{zw}|$;
scalar deviations are absolute differences.

\medskip\subsection{Intertwining and $\gamma$-independence}

For the biased random walk (Stratum~2, canonical) we compute,
for a set of $\gamma$ values, the coupling constant $C(\pi^*,\gamma)$
and the commutator norm $\|[\Delta^2\mathcal{R},R(\gamma)]\|$.
Table~\ref{tab:intertwining} presents the results for the
canonical configuration $a=10$, $p=0.6$, reset sites $\{3,7\}$.
In every case $C(\pi^*,\gamma)$ equals
$C^*=1/(1+\sqrt{K})=q_{a/2}^{(0)}=0.1163636364$ to machine precision
($|C-C^*|<10^{-16}$), independently of $\gamma$, and the commutator
norm is below $5\times10^{-14}$, confirming the resolvent
intertwining of Section~\ref{sec:intertwining} as an exact identity. The sum
$\sum_z\chi_z(\gamma)$ is zero to within $10^{-15}$, confirming
that $\chi(\gamma)$ lies in the tangent space of the simplex.

\begin{table}[t]
\centering
\caption{Intertwining and $\gamma$-independence for the biased
random walk ($a=10$, $p=0.6$, reset sites $\{3,7\}$).
$C^*=0.1163636364$, $K=(p/q)^a=57.66503906$, $\sqrt{K}=7.59375$.}
\begin{tabular}{cccc}
\toprule
$\gamma$ & $C(\pi^*,\gamma)$ & $\|[\Delta^2\mathcal{R},R(\gamma)]\|$ & $\sum_z\chi_z$ \\
\midrule
0.05 & 0.1163636364 & $3.77\times10^{-14}$ & $-3.55\times10^{-15}$ \\
0.15 & 0.1163636364 & $8.85\times10^{-15}$ & $0.00\times10^{0}$ \\
0.25 & 0.1163636364 & $9.76\times10^{-15}$ & $-3.55\times10^{-15}$ \\
0.35 & 0.1163636364 & $6.73\times10^{-15}$ & $0.00\times10^{0}$ \\
0.45 & 0.1163636364 & $9.99\times10^{-15}$ & $-8.88\times10^{-16}$ \\
0.55 & 0.1163636364 & $4.31\times10^{-15}$ & $-2.66\times10^{-15}$ \\
0.65 & 0.1163636364 & $8.66\times10^{-15}$ & $1.78\times10^{-15}$ \\
0.75 & 0.1163636364 & $1.01\times10^{-15}$ & $0.00\times10^{0}$ \\
0.85 & 0.1163636364 & $4.79\times10^{-16}$ & $0.00\times10^{0}$ \\
0.95 & 0.1163636364 & $2.82\times10^{-16}$ & $0.00\times10^{0}$ \\
\bottomrule
\end{tabular}
\label{tab:intertwining}
\end{table}

\medskip\subsection{Universality of $C^*$}

Table~\ref{tab:universality} verifies that $C^*$ depends only on
the invariant $K$, not on the choice of reset sites nor on the
specific chain. For fixed $(a,p)$, all symmetric site pairs give
identical $C^*$. Moreover, the non-homogeneous chain of
Example~3 (Stratum~1), with the same $K=(p/q)^a$ as the biased
random walk, yields the same $C^*$ --- confirming that $C^*$ is a
function of $K$ alone.

\begin{table}[t]
\centering
\caption{Universality of $C^*$ across symmetric site pairs and
across different chains. Top: invariance under site choice
($a=10$, $p=0.6$). Middle: $a=12$, $p=0.65$. Bottom: the
Stratum~1 chain (non-homogeneous, $a=6$, $p=0.6$) sharing the
same $C^*$ determined by its $K$.}
\begin{tabular}{cccc}
\toprule
$a$ & $p$ & Sites & $C(\pi^*,0.5)$ \\
\midrule
10 & 0.60 & $\{1,9\}$ & 0.1163636364 \\
10 & 0.60 & $\{2,8\}$ & 0.1163636364 \\
10 & 0.60 & $\{3,7\}$ & 0.1163636364 \\
10 & 0.60 & $\{4,6\}$ & 0.1163636364 \\
\midrule
12 & 0.65 & $\{2,10\}$ & 0.0237941145 \\
12 & 0.65 & $\{3,9\}$ & 0.0237941145 \\
12 & 0.65 & $\{4,8\}$ & 0.0237941145 \\
12 & 0.65 & $\{5,7\}$ & 0.0237941145 \\
\midrule
\multicolumn{4}{l}{\textit{Stratum~1, non-homogeneous ($a=6$, $p=0.6$):}} \\
6 & 0.60 & $\{1,5\}$ & 0.2285714286 \\
\bottomrule
\end{tabular}
\label{tab:universality}
\end{table}

The Stratum~1 chain has site-dependent rates ($p_1=0.30$,
$p_2=0.375$, $p_3=0.42$, $q_1=0.20$, $q_2=0.25$, $q_3=0.28$,
symmetric under $\sigma$), yet its coupling constant equals
$1/(1+\sqrt{(p/q)^6})=0.2285714286$ exactly, illustrating the
universality corollary of Sec.~\ref{sec:universality}.

\begin{figure*}[t]
\centering
\begin{tikzpicture}
\begin{axis}[
    name=panelA,
    width=8.4cm, height=6.2cm,
    xlabel={$\gamma$}, ylabel={$C(\pi,\gamma)$},
    xmin=0, xmax=1, ymin=0, ymax=1,
    xtick={0,0.25,0.5,0.75,1},
    ytick={0,0.2,0.4,0.6,0.8,1},
    grid=major, grid style={gray!25,line width=0.4pt},
    legend pos=north west, legend style={font=\scriptsize,fill=white,fill opacity=0.85,draw=gray!40},
    title={\small (a) neutral vs.\ non-neutral ($a{=}10$, $p{=}0.6$)},
    title style={yshift=-1pt},
]
\addplot[blue!70,thick,mark=none] coordinates {(0.020,0.33164)(0.080,0.47234)(0.140,0.60061)(0.200,0.70926)(0.260,0.79590)(0.320,0.86161)(0.380,0.90942)(0.440,0.94293)(0.500,0.96560)(0.560,0.98037)(0.620,0.98955)(0.680,0.99495)(0.740,0.99787)(0.800,0.99927)(0.860,0.99983)(0.920,0.99998)(0.980,1.00000)(0.980,1.00000)};
\addlegendentry{$\pi=\delta_3$}
\addplot[orange!85!black,thick,mark=none] coordinates {(0.020,0.26475)(0.080,0.34138)(0.140,0.40093)(0.200,0.44520)(0.260,0.47708)(0.320,0.49950)(0.380,0.51494)(0.440,0.52536)(0.500,0.53222)(0.560,0.53662)(0.620,0.53932)(0.680,0.54090)(0.740,0.54175)(0.800,0.54216)(0.860,0.54232)(0.920,0.54237)(0.980,0.54237)(0.980,0.54237)};
\addlegendentry{$\pi=0.8\delta_3{+}0.2\delta_7$}
\addplot[teal,thick,mark=none] coordinates {(0.020,0.03377)(0.080,0.01900)(0.140,0.01140)(0.200,0.00706)(0.260,0.00443)(0.320,0.00278)(0.380,0.00172)(0.440,0.00105)(0.500,0.00062)(0.560,0.00035)(0.620,0.00018)(0.680,0.00009)(0.740,0.00004)(0.800,0.00001)(0.860,0.00000)(0.920,0.00000)(0.980,0.00000)(0.980,0.00000)};
\addlegendentry{$\pi=\delta_7$}
\addplot[red,very thick] coordinates {(0.020,0.11636)(0.080,0.11636)(0.140,0.11636)(0.200,0.11636)(0.260,0.11636)(0.320,0.11636)(0.380,0.11636)(0.440,0.11636)(0.500,0.11636)(0.560,0.11636)(0.620,0.11636)(0.680,0.11636)(0.740,0.11636)(0.800,0.11636)(0.860,0.11636)(0.920,0.11636)(0.980,0.11636)(0.980,0.11636)};
\addlegendentry{$\pi=\pi^*$ (neutral)}
\draw[red,dashed,line width=0.6pt] (axis cs:0,0.11636) -- (axis cs:1,0.11636);
\node[red,font=\scriptsize,anchor=south east] at (axis cs:0.98,0.11636) {$C^*$};
\end{axis}
\begin{axis}[
    name=panelB, at={(panelA.east)}, anchor=west, xshift=1.1cm,
    width=8.4cm, height=6.2cm,
    xlabel={$\gamma$}, ylabel={$C(\pi^*,\gamma)$},
    xmin=0, xmax=1, ymin=0.1156, ymax=0.1172,
    xtick={0,0.25,0.5,0.75,1},
    ytick={0.1160,0.1164,0.1168},
    yticklabel style={/pgf/number format/fixed,/pgf/number format/precision=4},
    scaled y ticks=false,
    grid=major, grid style={gray!25,line width=0.4pt},
    legend pos=south east, legend columns=2, legend style={font=\scriptsize,fill=white,fill opacity=0.9,draw=gray!40,/tikz/every even column/.append style={column sep=4pt}},
    title={\small (b) universality: distinct chains, same $K$},
    title style={yshift=-1pt},
]
\addplot[red,only marks,mark=*,mark size=1.3pt] coordinates {(0.020,0.11636)(0.100,0.11636)(0.180,0.11636)(0.260,0.11636)(0.340,0.11636)(0.420,0.11636)(0.500,0.11636)(0.580,0.11636)(0.660,0.11636)(0.740,0.11636)(0.820,0.11636)(0.900,0.11636)(0.980,0.11636)(0.980,0.11636)};
\addlegendentry{homogeneous}
\addplot[blue!70,only marks,mark=square,mark size=1.3pt] coordinates {(0.020,0.11636)(0.100,0.11636)(0.180,0.11636)(0.260,0.11636)(0.340,0.11636)(0.420,0.11636)(0.500,0.11636)(0.580,0.11636)(0.660,0.11636)(0.740,0.11636)(0.820,0.11636)(0.900,0.11636)(0.980,0.11636)(0.980,0.11636)};
\addlegendentry{non-homog.\ 1}
\addplot[teal,only marks,mark=triangle,mark size=1.6pt] coordinates {(0.020,0.11636)(0.100,0.11636)(0.180,0.11636)(0.260,0.11636)(0.340,0.11636)(0.420,0.11636)(0.500,0.11636)(0.580,0.11636)(0.660,0.11636)(0.740,0.11636)(0.820,0.11636)(0.900,0.11636)(0.980,0.11636)(0.980,0.11636)};
\addlegendentry{non-homog.\ 2}
\addplot[orange!85!black,only marks,mark=diamond,mark size=1.6pt] coordinates {(0.020,0.11636)(0.100,0.11636)(0.180,0.11636)(0.260,0.11636)(0.340,0.11636)(0.420,0.11636)(0.500,0.11636)(0.580,0.11636)(0.660,0.11636)(0.740,0.11636)(0.820,0.11636)(0.900,0.11636)(0.980,0.11636)(0.980,0.11636)};
\addlegendentry{homog., sites $\{1,9\}$}
\draw[black,dashed,line width=0.7pt] (axis cs:0,0.116364) -- (axis cs:1,0.116364);
\node[black,font=\scriptsize,anchor=south west] at (axis cs:0.03,0.116364) {$C^*=1/(1+\sqrt{K})$};
\end{axis}
\end{tikzpicture}
\caption{Universality of the critical value $C^*$. \textbf{(a)} For a
fixed chain (biased random walk, $a=10$, $p=0.6$), the coupling
functional $C(\pi,\gamma)$ depends strongly on $\gamma$ for generic
reset distributions $\pi$ (blue, orange, teal), but is exactly
constant, $C(\pi^*,\gamma)\equiv C^*$, for the reset-neutral
distribution $\pi^*$ (red). The generic distributions
($\delta_3$, $0.8\,\delta_3+0.2\,\delta_7$, $\delta_7$) are chosen as
representative non-neutral examples, away from the separatrix, to
display the typical strong $\gamma$-dependence. \textbf{(b)} The neutral value is
\emph{universal}: four dynamically distinct $(\sigma,\kappa)$-reversible
chains sharing the same invariant $K=(p/q)^a$ --- one homogeneous, two
genuinely non-homogeneous (site-dependent rates obeying
$p_z=(p/q)q_{a-z}$), and the homogeneous chain with a different reset
pair --- all collapse onto the same $C^*=1/(1+\sqrt{K})$ (dashed) to
machine precision, independently of $\gamma$. The vertical scale spans
only $1.2\times10^{-3}$ around $C^*$.}
\label{fig:universality_collapse}
\end{figure*}
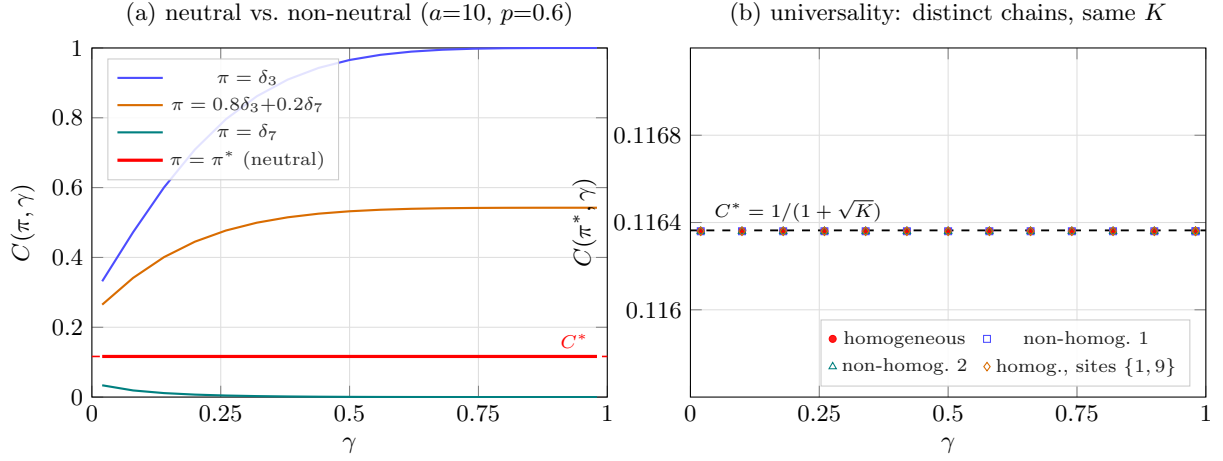

Figure~\ref{fig:universality_collapse} makes the universality visible.
Panel~(a) contrasts the strong $\gamma$-dependence of $C(\pi,\gamma)$
for generic reset distributions with the exact flatness of the
neutral curve; panel~(b) shows four dynamically distinct chains
sharing the same $K$ collapsing onto the single value $C^*$ to machine
precision.

\medskip\subsection{The antisymmetric sector: the representative $\chi_z(\gamma)$}

\begin{figure*}[t]
\centering
\begin{tikzpicture}
\begin{axis}[
    width=10cm,
    height=7cm,
    xlabel={$z$},
    ylabel={$\chi_z(\gamma)$},
    xmin=0.5, xmax=9.5,
    ymin=-14, ymax=14,
    xtick={1,2,3,4,5,6,7,8,9},
    ytick={-12,-8,-4,0,4,8,12},
    grid=major,
    grid style={line width=0.5pt, draw=gray!40},
    legend pos=north east,
    legend style={font=\small},
    title={Sector representative $\chi_z(\gamma)$ for $a=10$, $p=0.6$},
]
\addplot[thick, blue, mark=*, mark size=1.8pt] coordinates {
(1,11.8616)(2,12.5804)(3,9.3574)(4,4.8535)(5,0.0000)(6,-4.8535)(7,-9.3574)(8,-12.5804)(9,-11.8616)
};
\addlegendentry{$\gamma=0.2$}
\addplot[thick, red, mark=square*, mark size=1.8pt] coordinates {
(1,7.2010)(2,6.5053)(3,4.4207)(4,2.1944)(5,0.0000)(6,-2.1944)(7,-4.4207)(8,-6.5053)(9,-7.2010)
};
\addlegendentry{$\gamma=0.5$}
\addplot[thick, green!50!black, mark=triangle*, mark size=1.8pt] coordinates {
(1,5.5759)(2,4.3874)(3,2.8533)(4,1.3990)(5,0.0000)(6,-1.3990)(7,-2.8533)(8,-4.3874)(9,-5.5759)
};
\addlegendentry{$\gamma=0.8$}
\draw[dashed, gray!60] (axis cs:5,-14) -- (axis cs:5,14);
\draw[dashed, gray!60] (axis cs:0.5,0) -- (axis cs:9.5,0);
\end{axis}
\end{tikzpicture}
\caption{Antisymmetric sector representative $\chi_z(\gamma)$ versus site $z$ for
three values of $\gamma$ ($a=10$, $p=0.6$). The antisymmetry
$\chi_{a-z}=-\chi_z$ and the exact node at the midpoint $z=5$ are
visible. The amplitude decreases as $\gamma$ increases. Dashed
lines mark the midpoint and the zero axis.}
\label{fig:psi_field}
\end{figure*}
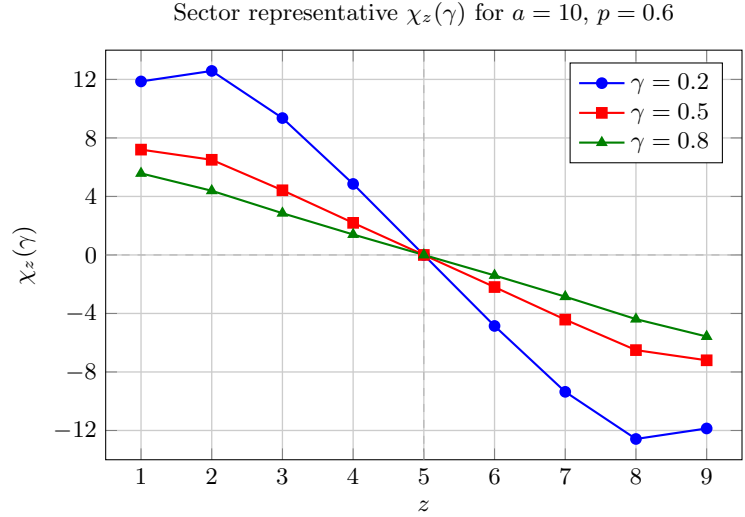

Figure~\ref{fig:psi_field} shows the canonical sector
representative
\[
\chi_z(\gamma) = [R(\gamma)\xi]_z - [R(\gamma)\xi]_{\sigma(z)},
\qquad \xi_z = \sqrt{\kappa(z)},
\]
for the biased random walk with $a=10$, $p=0.6$, at three values
of $\gamma$. Here $\xi_z=\sqrt{\kappa(z)}=(p/q)^{(a-z)/2}$ is the
square root of the canonical weight, so $\chi=(I-\mathcal{R})R\xi$ with
$(\mathcal{R}f)_z=f_{\sigma(z)}$ the reflection operator. The field is
exactly $\mathcal{R}$-antisymmetric ($\chi_{\sigma(z)}=-\chi_z$),
and its amplitude decreases monotonically as $\gamma$ increases,
reflecting the suppression of the orientation under stronger
resetting. The orientation field proper,
$\psi(\gamma)=R(\gamma)(b^{(0)}-C^*b)$, shares the exact midpoint
node and satisfies the $\Delta$-weighted antisymmetry
$\Delta^{-1}\psi(\sigma(z))=-\Delta^{-1}\psi(z)$; we display
$\chi$ because it is the plainly antisymmetric representative of
the same sector.

\medskip
A structural feature deserves emphasis. When $a$ is even, the
midpoint $z=a/2$ is the unique fixed point of the involution
($\sigma(a/2)=a/2$), and therefore
\[
\begin{aligned}
\chi_{a/2}(\gamma)
&= [R(\gamma)\xi]_{a/2}-[R(\gamma)\xi]_{\sigma(a/2)}=0,\\[2pt]
&\quad\text{and likewise }\psi_{a/2}(\gamma)=0,
\qquad \forall\,\gamma\in(0,1).
\end{aligned}
\]
This is an \emph{exact} identity, not a numerical approximation:
the midpoint is a genuine antisymmetric node of the resolvent
geometry, independent of $\gamma$, of $\xi$, and of the chain.
It identifies the separatrix as a structural node --- the
informational equator of Sec.~\ref{sec:universality}. (For $a$ odd the involution
has no integer fixed point and the two central sites are exchanged
instead.)

Table~\ref{tab:psi_values} provides the numerical values of
$\chi_z(\gamma)$ together with the reference vector
$\xi_z=\sqrt{\kappa(z)}$. The antisymmetry is exact to the
reported precision.

\begin{table*}[t]
\centering
\caption{Numerical values of $\chi_z(\gamma)$ for the biased
random walk ($a=10$, $p=0.6$). The vector $\xi_z=\sqrt{\kappa(z)}$
is independent of $\gamma$.}
\begin{tabular}{ccccc}
\toprule
$z$ & $\chi_z(0.2)$ & $\chi_z(0.5)$ & $\chi_z(0.8)$ & $\xi_z$ \\
\midrule
1 & 11.86162445 & 7.20103516 & 5.57590342 & 6.20027091 \\
2 & 12.58037209 & 6.50533436 & 4.38737318 & 5.06250000 \\
3 & 9.35743057 & 4.42074625 & 2.85330335 & 4.13351394 \\
4 & 4.85354334 & 2.19438312 & 1.39904731 & 3.37500000 \\
5 & 0.00000000 & 0.00000000 & 0.00000000 & 2.75567596 \\
6 & -4.85354334 & -2.19438312 & -1.39904731 & 2.25000000 \\
7 & -9.35743057 & -4.42074625 & -2.85330335 & 1.83711731 \\
8 & -12.58037209 & -6.50533436 & -4.38737318 & 1.50000000 \\
9 & -11.86162445 & -7.20103516 & -5.57590342 & 1.22474487 \\
\bottomrule
\end{tabular}
\label{tab:psi_values}
\end{table*}

\medskip\subsection{Summary}

The numerical experiments confirm, to machine precision, the
three central claims of the paper: (i) the resolvent intertwining
$[\Delta^2\mathcal{R},R(\gamma)]=0$ holds exactly, and hence
$C(\pi^*,\gamma)=C^*$ is independent of $\gamma$; (ii) the critical
value $C^*$ depends only on the invariant $K$, not on the reset
sites nor on the chain within the $(\sigma,\kappa)$-reversible
class; and (iii) the sector representative $\chi(\gamma)$ exhibits the
predicted antisymmetric structure with an exact node at the
midpoint, its amplitude decreasing with $\gamma$. Importantly, no
parameter fitting is involved: every quantity --- $K$, $C^*$, the
operator $\Delta^2\mathcal{R}$, and the fields $\psi(\gamma)$, $\chi(\gamma)$ --- is
determined analytically from $(p,q,a)$ and the symmetry structure,
and the experiments verify exact identities rather than calibrate
free constants. These results
provide strong numerical confirmation of the framework of
Sections~\ref{sec:resolvent}--\ref{sec:universality}.

\section{Conclusions and perspectives}
\label{sec:conclusions}

In this paper we have uncovered the resolvent origin of the spectral
duality that governs reset-neutral distributions in Markov chains
with geometric resetting. Starting from the abstract conditions of
Paper~III, we showed that the spectral duality
$B_\nu(z)=\kappa(z)\,A_\nu(\sigma(z))$ is equivalent to an
intertwining relation for the resolvent
$R(\gamma)=(I-(1-\gamma)P)^{-1}$,
\[
[\Delta^2\mathcal{R},\,R(\gamma)]=0,
\]
where $\mathcal{R}$ is the reflection operator of the involution
$\sigma$ and $\Delta=\operatorname{diag}(\sqrt{\kappa(z)})$. This
resolvent symmetry underlies both the existence of reset-neutral
distributions and the universal critical value
$C^*=1/(1+\sqrt{K})$, with $K=\kappa(z)\kappa(\sigma(z))$.

\medskip\noindent
We emphasize what we regard as the central conceptual contribution
of this paper. Papers~I--III established the existence of
reset-neutral distributions, the invariant $C^*$, and the geometry
of the separatrix. The present work identifies their common origin:
reset-neutrality is the manifestation of a hidden symmetry of the
resolvent, $[\Delta^2\mathcal{R},R(\gamma)]=0$. This reframes a
phenomenon that appears specific to resetting --- the
$\gamma$-independence of $C(\pi^*,\gamma)$ --- as a structural
property of the resolvent that holds for every reset rate
simultaneously. It is this shift, from a resetting phenomenon to a
resolvent symmetry, that distinguishes the present analysis from
Papers~I--III and opens the abstract $(\sigma,\kappa,K)$ framework
to chains far beyond the biased random walk.

\medskip\noindent\textbf{Summary of main results.}
\begin{enumerate}
\item \textbf{Resolvent intertwining as the core symmetry.}
(S3) is equivalent to $[\Delta^2\mathcal{R},R(\gamma)]=0$, and the
normalized operator $\widetilde{\mathcal{T}}=K^{-1/2}\Delta^2\mathcal{R}$
is an involution ($\widetilde{\mathcal{T}}^2=I$, since
$(\Delta^2\mathcal{R})^2=K\,I$ by (H3)). This gives a clean
operator-theoretic form of the twisted symmetry.
\item \textbf{The exact neutrality identity and the orientation field.}
The deviation of the coupling functional from its critical value
satisfies the exact identity
$(I+\mathcal{R})\,\Delta^{-1}R(\gamma)(b^{(0)}-C^*b)=0$ for every
$\gamma$, which forces $C^*=1/(1+\sqrt{K})$ mode by mode. The
orientation field is the deviation itself,
\[
\psi(\gamma)=u-C^*s=R(\gamma)\,(b^{(0)}-C^*b),
\]
with $\Delta^{-1}\psi(\gamma)$ $\sigma$-antisymmetric; it obeys the
exact sign law
$\operatorname{sgn}(C(\pi,\gamma)-C^*)
=\operatorname{sgn}\langle\pi,\psi(\gamma)\rangle$ and is confined
to the fixed subspace $\Delta\cdot\ker(I+\mathcal{R})$ of dimension
$n_{\mathrm{pair}}$. The plainly antisymmetric representative
$\chi(\gamma)=(I-\mathcal{R})R(\gamma)\xi$, $\xi_z=\sqrt{\kappa(z)}$,
is the computational companion of the same sector.
\item \textbf{The universal class and its stratification.}
The $(\sigma,\kappa)$-reversible chains are the natural setting for
(S3). For nearest-neighbor chains with $\sigma(z)=a-z$, the rate
condition $p_z=(\kappa(z)/\kappa(z+1))\,q_{a-z}$ holds together with
a $\sigma$-symmetric holding profile, and the class stratifies into
Stratum~1 (lazy versions of the biased walk, same $K$ and $C^*$),
Stratum~2 (the biased random walk, canonical element of
Papers~I--III), and Stratum~3 (non-exponential weights, which occur
already for nearest-neighbor chains, and longer-range chains).
Because the duality holds step by step in time, the reset-neutral
distributions are neutral for every resetting mechanism, not only for
geometric resetting.
\item \textbf{Universality and $h$-transform duality.}
$C^*=1/(1+\sqrt{K})$ depends only on the scalar $K$ --- not on
$\gamma$, the reset distribution, or the chain within the class.
The Doob $h$-transform induces $K\mapsto1/K$ and hence
$C^*\mapsto1-C^*$, with fixed point $C^*=1/2$ realizing the
separatrix as an \emph{informational equator}. This duality is
itself an involution on parameter space ($K\mapsto1/K$ applied
twice restores $K$), complementing the two other involutions of
the theory: $\sigma$ on the state space and
$\widetilde{\mathcal{T}}=K^{-1/2}\Delta^2\mathcal{R}$ on
operators. The reset-neutral structure is thus organized by a
triad of involutions acting at the levels of states, operators,
and parameters.
\item \textbf{Numerical verification.}
All claims were confirmed to machine precision: the exact
$\gamma$-independence of $C(\pi^*,\gamma)=C^*$ ($<10^{-16}$); the
universality of $C^*$ across site pairs and across the Stratum~1
chain; the antisymmetric orientation field with an exact midpoint
node; and the vanishing commutator ($\sim10^{-14}$).
\end{enumerate}

\medskip\noindent\textbf{Conceptual reading.}
Geometric resetting is governed not by the resetting mechanism
alone but by a dual geometry encoded in the involution $\sigma$ and
the weight $\kappa$. The invariant $K$ acts as an orientation
parameter for the spectral geometry: $K<1$ and $K>1$ correspond to
opposite drift orientations, exchanged by the $h$-transform, with
$K=1$ the self-dual equilibrium. The universality of $C^*$ reflects
the rigidity of the underlying $(\sigma,\kappa)$ structure.

\medskip\noindent\textbf{Open questions and future work.}

\textit{Non-exponential weights.} For nearest-neighbor chains with
the reflection involution, the weight $\kappa$ need not be
exponential, but it must be monotone (Remark~\ref{rem:gauge} and
Appendix~\ref{app:nn}). How the spectral and response properties of
the chains of the class depend on the weight is left open.

\textit{Operator-theoretic foundations.} A deeper interpretation in
terms of the generator and its Perron--Frobenius eigenvector is
expected to yield a unified framework.

\textit{Information geometry (Paper~V~\cite{PaperV}).} The representation
$\psi(\gamma)=R(\gamma)(b^{(0)}-C^*b)$, together with the finite
reduction $\psi(\gamma)\in\Delta\cdot\ker(I+\mathcal{R})$ of
dimension $n_{\mathrm{pair}}$, is the starting point for
Paper~V, where the simplex carries the Fisher--Rao metric and the
separatrix $\Sigma$ is a totally geodesic submanifold of
codimension $n_{\mathrm{pair}}$. There the position of $\pi$ relative to
$\Sigma$ reduces to a signed distance exactly when
$n_{\mathrm{pair}}=1$, the case in which $\Sigma$ separates the
simplex into two sides. That separation is necessary but not
sufficient for a global orientation law, and two further conditions
enter. The first is a scalar sign condition: the Wronskians
$\mathcal{W}[\varphi,s(z;\cdot)]$ must share a common strict sign
$\varepsilon$ across the reset sites, neutral sites included. Under
$n_{\mathrm{pair}}=1$ and that condition, Paper~V establishes a
sign law for $\partial_\gamma C$ referred to the fixed director
$\mathbf{n}$ and the constant $\varepsilon$. The second condition is
a calibration: writing $\psi(\gamma)=\varphi(\gamma)\mathbf{n}$, the
sign law just described coincides with the identification considered
in Paper~III --- which is referred to $\psi(\gamma)$ itself --- only
when $\varepsilon=\operatorname{sgn}\varphi$. That calibration is not
implied by the structural conditions: there exist structures with
$n_{\mathrm{pair}}=1$ satisfying the sign condition for which the
identification nevertheless fails, the two orientations
being opposite. The canonical family satisfies the calibration, so
for it the identification is a theorem; in the abstract class it is
an additional requirement, and Paper~V isolates it as such.

\textit{Applications.} The explicit characterization of reset-neutral
distributions and the stratification of the class enable the design
of non-homogeneous chains sharing the reset-invariant properties of
the biased random walk, with potential use in optimal resetting,
controlled random walks, and queueing
systems~\cite{Reuveni2016,PalReuveni2017,Asmussen2003}.

\medskip\noindent\textbf{Closing remark.}
The resolvent intertwining $[\Delta^2\mathcal{R},R(\gamma)]=0$
reveals that the reset-neutral phenomenon is not an accident of the
biased random walk but a structural consequence of a hidden
symmetry of the resolvent. The $(\sigma,\kappa)$-reversible class
provides a rich family of processes sharing the same universal
$C^*$, and the orientation field $\psi(\gamma)$, now explicitly
represented through the resolvent, is the bridge to the
information-geometric theory of Paper~V.

\appendix
\section{Nearest-neighbour realizations of the class}
\label{app:nn}

We consider nearest-neighbour chains on $\mathcal{T}=\{1,\dots,a-1\}$
with up-rates $p_z>0$, down-rates $q_z>0$ and holding probabilities
$d_z=1-p_z-q_z\ge0$, the reflection $\sigma(z)=a-z$, and a weight
$\kappa$ normalized by the twisted boundary pairing (3.12). It is
convenient to set $\kappa(0):=K$ and $\kappa(a):=1$, which is
consistent with $\kappa(z)\kappa(a-z)=K$, and to write
$r_z:=\kappa(z)/\kappa(z+1)$ for $z=0,\dots,a-1$. The constant-product
relation gives $r_{a-z}=r_{z-1}$ for $z=1,\dots,a$, and
$K=\kappa(0)/\kappa(a)=\prod_{z=0}^{a-1}r_z$. By
Section~\ref{sec:universality}, $(\sigma,\kappa)$-reversibility together
with (3.12) reads
\begin{equation}
p_z=r_z\,q_{a-z}\quad(z=1,\dots,a-1),\qquad d_z=d_{a-z},
\tag{A1}
\end{equation}
the case $z=a-1$ being the boundary pairing
$p_{a-1}=\kappa(a-1)\,q_1$.

As shown by T.~Newton~\cite{Newton2026}, these conditions do not force
the weight to be exponential: for every even $a\ge6$ and every
strictly decreasing weight with constant product $K$ there are chains
realizing it (Remark~\ref{rem:gauge}). The next proposition shows that
monotonicity is also necessary: no other weights can occur, and for
even $a\ge6$ and $K>1$ the weights that occur are exactly the strictly
decreasing ones.

\medskip\noindent\textbf{Proposition A.1 (Monotonicity of the
weight).} \textit{If a nearest-neighbour chain with positive rates
satisfies (A1), then $\kappa$ is strictly decreasing on
$\{0,\dots,a\}$ if $K>1$, strictly increasing if $K<1$, and constant
if $K=1$.}

\textit{Proof.} By (A1) and $r_{a-z}=r_{z-1}$, $p_z=r_zq_{a-z}$ and
$p_{a-z}=r_{z-1}q_z$. The symmetry of the holding probabilities,
$p_z+q_z=p_{a-z}+q_{a-z}$, then becomes
\[
q_{a-z}\,(r_z-1)=q_z\,(r_{z-1}-1),\qquad z=1,\dots,a-1 .
\]
Since $q_z,q_{a-z}>0$, the numbers $r_z-1$ and $r_{z-1}-1$ are both
positive, both negative, or both zero, for every $z=1,\dots,a-1$.
Hence all of $r_0,\dots,r_{a-1}$ lie strictly on the same side of $1$,
or all equal $1$. In the first case $\kappa$ is strictly monotone on
$\{0,\dots,a\}$, decreasing if the $r_z$ exceed $1$ and increasing
otherwise, and $K=\prod_z r_z$ is accordingly larger or smaller than
$1$; in the second case $\kappa$ is constant and $K=1$.
\hfill$\square$

\medskip
In particular, the exponential weight $\kappa(z)=(p/q)^{a-z}$ of the
biased random walk is the case of constant ratios $r_z\equiv p/q$, and
no non-monotone weight is realized by a nearest-neighbour chain of
the class.

The critical value admits a description that does not involve
resetting at all. Let $h(z)$ denote the classical ruin probability of
the chain without resetting, the probability of absorption at
$\partial_0$ before $\partial_1$ from $z$.

\medskip\noindent\textbf{Proposition A.2 (Reset-free meaning of
$C^*$).} \textit{Let $P$ be $(\sigma,\kappa)$-reversible and satisfy
(3.12).}
\begin{itemize}
\item[(i)] \textit{For every reset-neutral distribution $\pi^*$,
$\sum_z\pi^*_z\,h(z)=C^*=1/(1+\sqrt{K})$: the critical value is the
ruin probability of the chain without resetting started from
$\pi^*$. If $\sigma$ has a fixed point $z_0$ (for $\sigma(z)=a-z$, the
midpoint $a/2$ when $a$ is even), then $h(z_0)=C^*$.}
\item[(ii)] \textit{For nearest-neighbour chains with
$\sigma(z)=a-z$,}
\[
K=\frac{p_{a-1}}{q_1}\prod_{z=1}^{a-1}\frac{p_z}{q_z},
\qquad
C^*=\frac{1}{1+\sqrt{K}} .
\]
\end{itemize}

\textit{Proof.} (i) By the Remark on the independence of the resetting
mechanism (Section~\ref{sec:universality}), a chain started from
$\pi^*$ and reset to $\pi^*$ has ruin probability $C^*$ under every
resetting mechanism under which absorption is almost sure. The empty
mechanism, with no resets, is one of them, since absorption is certain
from every transient state; its ruin probability is
$\sum_z\pi^*_zh(z)$. If $\sigma(z_0)=z_0$, then $\kappa(z_0)=\sqrt K$
and $\pi^*=\delta_{z_0}$ is reset-neutral, which gives $h(z_0)=C^*$.
(ii) By (A1), $r_z=p_z/q_{a-z}$ for $z=1,\dots,a-1$, and
$r_0=r_{a-1}=\kappa(a-1)=p_{a-1}/q_1$. Hence
$K=\prod_{z=0}^{a-1}r_z
=(p_{a-1}/q_1)\prod_{z=1}^{a-1}p_z\big/\prod_{z=1}^{a-1}q_{a-z}$,
and the last product equals $\prod_{z=1}^{a-1}q_z$.
\hfill$\square$

\medskip
For the biased random walk, (i) is the identity
$C^*=q^{(0)}_{a/2}$ of Papers~I and~II and (ii) gives $K=(p/q)^a$.
The proposition shows that both are properties of the whole class:
for every nearest-neighbour chain of the class with $a$ even, whether
or not its weight is exponential, the reset-neutral value is the
classical ruin probability from the midpoint, and it is determined by
the product of the rate ratios. For the chain of
Remark~\ref{rem:gauge}, $h(3)=1/5=1/(1+\sqrt{16})$, and the product
formula gives $K=(0.2/0.1)\cdot8=16$.

\begin{acknowledgments}
The author thanks T.~Newton for correspondence in which explicit
constructions satisfying the structural assumptions of Paper~III
were shown to violate the global orientation of the response in the
form considered there. Those examples make clear that the sign
condition entering the statements above is genuinely necessary rather
than merely convenient, and prompted the corresponding sharpening of
the forward-looking remarks in Sections~\ref{sec:intro}
and~\ref{sec:conclusions}. He also thanks T.~Newton for showing that
the structural conditions do not force the weight $\kappa$ to be
exponential, through the family of examples recorded in
Remark~\ref{rem:gauge}.
\end{acknowledgments}

\end{document}